\documentclass{amsart}

\newtheorem{theorem}{Theorem}[section]
\newtheorem{lemma}[theorem]{Lemma}
\newtheorem{proposition}[theorem]{Proposition}
\newtheorem{corollary}[theorem]{Corollary}
\newtheorem{theoremi}{Theorem}

\theoremstyle{definition}

\newtheorem{remark}[theorem]{Remark}
\newtheorem*{acknowledgement}{Acknowledgement}

\theoremstyle{remark}

\newcommand\mynote[1]{\marginpar{\ \\ \small \tt #1}}
\newcommand\bel[1]{{\mynote{#1}}\begin{equation}\label{#1}}
\newcommand\mylabel[1]{\label{#1}}

\newcommand{\ZZ}{\mathbb{Z}}
\newcommand{\QQ}{\mathbb{Q}}

\newcommand{\FF}{\mathbb{F}}

\newcommand{\PP}{\mathbb{P}}
\renewcommand{\AA}{\mathbb{A}}

\newcommand  {\shF}     {\mathcal{F}}

\newcommand  {\shHom}   {\mathcal{H}\!\text{\textit{om}}}
\newcommand  {\shI}     {\mathcal{I}}
\newcommand  {\shJ}     {\mathcal{J}}

\newcommand  {\shL}     {\mathcal{L}}

\newcommand  {\foC}     {\mathfrak{C}}

\newcommand  {\foX}     {\mathfrak{X}}

\newcommand  {\foZ}     {\mathfrak{Z}}

\newcommand  {\Aut}     {\operatorname{Aut}}

\newcommand  {\Cl}      {\operatorname{Cl}}

\newcommand  {\coker}   {\operatorname{coker}}

\newcommand  {\disc}    {\operatorname{disc}}
\newcommand  {\Div}     {\operatorname{Div}}

\newcommand  {\edim}    {\operatorname{edim}}

\newcommand  {\Gal}     {\operatorname{Gal}}

\newcommand  {\Hom}     {\operatorname{Hom}}

\newcommand  {\id}      {\operatorname{id}}

\renewcommand  {\k}     {\kappa}
\renewcommand  {\ker }  {\operatorname{kern}}

\newcommand  {\Lie}     {\operatorname{Lie}}

\newcommand  {\lra}     {\longrightarrow}

\newcommand  {\maxid}   {\mathfrak{m}}

\newcommand  {\mult}    {\operatorname{mult}}

\newcommand  {\NS}      {\operatorname{NS}}

\renewcommand{\O}       {\mathcal{O}}

\newcommand  {\Pic}     {\operatorname{Pic}}

\newcommand  {\pr}      {\operatorname{pr}}
\newcommand  {\Proj}    {\operatorname{Proj}}

\newcommand  {\quadand} {\quad\text{and}\quad}

\newcommand  {\ra}      {\rightarrow}

\newcommand  {\red}     {{\operatorname{red}}}

\newcommand  {\sep}     {{\operatorname{sep}}}

\newcommand  {\Sing}    {\operatorname{Sing}}

\newcommand  {\Spec}    {\operatorname{Spec}}

\newcommand  {\Tor}     {\operatorname{Tor}}

\def\mydate{\number\day\space\ifcase\month \or January\or February\or March\or 
April\or May\or June\or July\or
August\or September\or October\or November\or December\fi \space\number\year}

\newcommand{\pd}  {\operatorname{pd}}

\newcommand{\lieg}{\mathfrak{g}}
\newcommand{\lieh}{\mathfrak{h}}
\newcommand{\lieD}{\mathfrak{D}}
\newcommand{\liea}{\mathfrak{a}}
\newcommand{\lieb}{\mathfrak{b}}
\newcommand{\liegl}{\mathfrak{gl}}
\newcommand{\End} {\operatorname{End}}
\newcommand{\Der} {\operatorname{Der}}
\newcommand{\Diff} {\operatorname{Diff}}
\newcommand{\homdim}{\operatorname{hom-dim}}
\newcommand{\Tordim}{\operatorname{Tor-dim}}

\newcommand{\flatdim} {\operatorname{flat-dim}}
\newcommand{\Res} {\operatorname{Res}}

\usepackage{amscd,amssymb,graphics}

\begin{document}

\title[Kummer surfaces]
      {Kummer surfaces for the selfproduct of the cuspidal rational curve}

\author[Stefan Schroer]{Stefan Schr\"oer}
\address{Mathematisches Institut, Heinrich-Heine-Universit\"at,
40225 D\"usseldorf, Germany}
\curraddr{}
\email{schroeer@math.uni-duesseldorf.de}


\dedicatory{Revised version, 11 November 2005}

\begin{abstract}
The classical Kummer construction attaches to an abelian surface
a K3 surface. As Shioda and Katsura showed, this construction breaks
down for supersingular abelian surfaces in characteristic two.
Replacing supersingular abelian surfaces by the selfproduct of
the rational cuspidal curve, and the sign involution by suitable
infinitesimal group scheme actions, I give   the correct  Kummer-type construction for
this situation. We encounter rational double points
of type $D_4$ and $D_8$, instead of type $A_1$.
It turns out that the resulting surfaces are supersingular K3 surfaces with Artin invariant one and two.
They lie in a 1-dimensional family obtained by simultaneous resolution, which
exists  after purely inseparable base change.
\end{abstract}

\maketitle
\tableofcontents

\section*{Introduction}

Let $A$ be an abelian surface over the complex numbers,
and $\iota:A\ra A$ the sign involution. The quotient surface
$Z=A/\iota$ is a normal surface with 16 rational double points of type $A_1$,
whose minimal resolution $X$ is a K3 surface. One also says that
$X$ is a Kummer K3 surface;  they play a fairly central role
in the theory of all complex K3 surfaces.

It is easy to see that the Kummer construction works in positive
characteristics $p\neq 2$ as well.
In contrast, Shioda \cite{Shioda 1974}  and Katsura  \cite{Katsura 1978} observed that the Kummer construction
breaks down in characteristic $p=2$ for supersingular abelian surfaces $A$.
In this case, they showed that
singularities on the quotient surface $Z$ are   elliptic singularities, and that the minimal
resolution $X$ is a rational surface.

The goal of this paper is to give a new type of Kummer construction
in the supersingular situation at $p=2$.
To explain this construction, let me   discuss the case
where $A$ is superspecial, that is, isomorphic to the selfproduct
$E\times E$ of supersingular elliptic curves.
My idea   is to replace the supersingular elliptic curve $E$
by a cuspidal rational curve $C$, and the group action of $\ZZ/2\ZZ$
by a suitable group scheme action of the infinitesimal group scheme $\alpha_2$.
In particular, we start with the \emph{nonnormal} surface $Y=C\times C$.
Nevertheless, it turns out that the quotient $Z=Y/\alpha_2$ is a normal surface,
whose singularities are rational double points of type $D_4$, $B_3$, or $D_8$,
at least if the ground field is perfect.
The first main result of this paper is the following:

\begin{theoremi}
The  minimal resolution of singularities $X\ra Z$ is a K3 surface.
\end{theoremi}

There are complications for nonperfect ground fields. It seems that
the minimal resolution could be a regular surface with trivial canonical class
that is not geometrically regular.
This  seems   to be an interesting subject matter in its own right, but
I do not pursue this topic in the present paper.

The Kummer construction also plays an important role in the theory of 
supersingular K3 surfaces in characteristic $p>0$.
Recall that a K3 surface $X$ is called \emph{supersingular} (in the sense of Shioda)
if its Picard number equals the second Betti number, that is, $\rho=22$.
Artin \cite{Artin 1962} introduced for such    K3 surfaces   an integer invariant
$1\leq\sigma_0\leq 10$ called the \emph{Artin invariant}, which can be defined
in terms of discriminants of the intersection form.
Shioda  \cite{Shioda 1978} proved for $p\neq 2$ that the supersingular K3 surfaces with Artin invariant
$1\leq\sigma_0\leq 2$ are precisely the Kummer K3 surfaces coming from  supersingular abelian
surfaces.  The second main result of this paper is:
 
 \begin{theoremi}
Our K3 surfaces $X$ are supersingular with Artin invariant $1\leq \sigma_0\leq 2$.
 \end{theoremi}
 
This depends on an analysis of the degenerate fibers in the quasielliptic fibrations
$f:X\ra\PP^1$ induced from the projections on $Y=C\times C$.
Drops in the Artin invariants are due to confluence of a pair of $D_4$-singularities
into a single $D_8$-singularity. 

Oort  \cite{Oort 1975} showed that any supersingular abelian surface is
an infinitesimal quotient of a superspecial abelian surface.
This implies that supersingular abelian surfaces $A$ form a 1-dimensional family,
in which the action of $\ZZ/2\ZZ$ is constant.
In our new Kummer construction, it is the other way round:
The nonnormal surface $Y=C\times C$ does not move, but the
action of $\alpha_2$ lies in a moving  family. This is reminiscent of the Moret-Bailly construction
\cite{Moret-Bailly 1981}.
Our construction gives a family $\foZ\ra S$ of  normal surfaces with rational double points over the punctured affine plane $S=\AA^2-0$. According to the work of Brieskorn \cite{Brieskorn 1968} and Artin \cite{Artin 1974b},
simultaneous resolutions exist rarely without base change. 
This is indeed the case in our situation:
 
\begin{theoremi}
The family of normal surface $\foZ\ra S$ admits a simultaneous minimal resolution of singularities
after certain purely inseparable base change.
\end{theoremi}
 
We achieve simultaneous resolution by successively blowing up Weil divisors in  multiple fibers on the
quasielliptic fibrations $Z\ra\PP^1$. The purely inseparable base change enters the picture,
because the original family contains integral fibers switching to multiple fibers over the perfect closure.
It turns out that the  resulting family of K3 surfaces $\foX\ra S'$, where $S'\ra S$ is the purely inseparable base change,
is induced from a family of $K3$ surfaces parameterized by the projective line $\PP^1$.
This is in accordance with results of   Rudakov and Shafarevich \cite{Rudakov; Safarevic 1978} and  Ogus \cite{Ogus 1983}
on the moduli of marked  supersingular K3 surfaces.

The paper is organized as follows:
Section \ref{generalities} discusses the  interplay between restricted Lie algebras
and infinitesimal group schemes of height $\leq 1$, in particular $\alpha_p$.
Section \ref{quotients} contains some results on the behavior of quotients
with respect to base change and singularities. In Section \ref{cuspidal} we examine
the rational cuspidal curve $C$ in characteristic two, determine its restricted Lie algebra, and read off
all possible $\alpha_2$-actions.
In Section \ref{selfproduct}, we consider the selfproduct $Y=C\times C$,
define   $\alpha_2$-actions on it, and establish some basic properties.
This is refined in the following two sections, where we determine
the structure of singularities on the  quotient $Z=Y/\alpha_2$.  
I return to global properties of $Y$ in Section \ref{quasielliptic},
where we analyse quasielliptic fibrations.
In Section \ref{k3} we turn to the  the resolution of singularities $X$,
and establish that it is a supersingular K3 surface provided that the parameters do not vanish
simultaneously. In Section \ref{discriminants}, we determine the Artin invariant.
Here we crucially use the structure of the singularities and the quasielliptic fibration.
Sections \ref{blowing}   contains some material on blowing up
rational surface singularities along curves. We use this in
Section \ref{genus}, where we consider blowing ups of Weil divisors
in genus-one fibrations. This is crucial in Section \ref{simultaneous}, where
we view  our normal K3 surfaces as lying in a family $\foZ\ra S$, and show that
this family admits a simultaneous resolution $\foX'\ra S'$ after a purely inseparable base change $S'\ra S$.
In the last section  we    show that this family is induced
from a family $\foX\ra\PP^1$.

\begin{acknowledgement}
I wish to thank Torsten Ekedahl and Nick Shepherd-Barron for helpful discussion. 
They pointed out to me the necessity of purely inseparable base change to obtain
simultaneous resolution of singularities.
I also thank Otto Kerner for helpful discussions on restricted Lie algebras,
Eckart Viehweg for discussions on moduli spaces, 
and Le Van Schr\"oer for Magma programming and careful proofreading.
Moreover, I   wish to thank the referee for writing a thorough report and for pointing out some mistakes.
\end{acknowledgement}

\section{Generalities on $\alpha_p$-actions}
\mylabel{generalities}

Let $k$ be a ground field of characteristic $p>0$. Then there is 
a  group scheme $\alpha_p$ whose values on $k$-algebras $R$ are
$\alpha_p(R)=\left\{r\in R\mid r^p=0\right\}$,
with addition as group law. As a scheme, we have
$\alpha_p=\Spec k[t]/(t^p)$, which is finite and infinitesimal.
Such group schemes exist only in characteristic $p>0$; it is 
frequently possible to explain    characteristic-$p$-phenomena in terms
of $\alpha_p$-actions.
In this section I collect some well-known facts on finite infinitesimal group schemes,
and in particular on $\alpha_p$. Proofs are   omitted.  
The  book of Demazure and Gabriel \cite{Demazure; Gabriel 1970} is an exhaustive reference.

The Lie algebra $\Lie(\alpha_p)=\alpha_p(k[\epsilon])$ is 
the 1-dimensional vector space $k\epsilon\subset k[\epsilon]$. The generator $\epsilon$
corresponds to the derivation $\partial:k[t]/(t^p)\ra k$ given by
the Kronecker delta $\partial(t^i)=\delta_{i1}$. 
The interpretation   $\Lie(\alpha_p)=\Der(k[t]/(t^p),k)$ shows that both
the Lie bracket
$[D,D']=D\circ D'-D'\circ D$ and the $p$-th power operation $D^{[p]}=D\circ\ldots\circ D$ 
($p$-fold composition) vanish.
In other words, $\Lie(\alpha_p)$ is a restricted Lie algebra
whose Lie bracket and $p$-map both vanish.

Recall that a \emph{restricted Lie algebra} over $k$ is a Lie algebra $\lieg$
endowed with a  $p$-map $\lieg\ra\lieg$,  $x\mapsto x^{[p]}$
(sometimes restricted Lie algebras are also called \emph{$p$-Lie algebras}).
By definition,  $p$-maps  satisfy three axioms. The first axiom relates them with scalar multiplication:
$$
(\lambda x)^{[p]} =\lambda^p x^{[p]} \quad\text{for all $\lambda\in k$, $x\in \lieg$}.
$$
The second axiom relates $p$-maps with Lie brackets:
$$
[x^{[p]},y]=[x,[x,[\ldots,[x,y]]]] \quad\text{for all $x,y\in\lieg$},
$$
where the right-hand side is a $p$-fold iterated Lie bracket. The third axiom relates
$p$-maps to addition:
$$
(x+y)^{[p]}=x^{[p]} + \Lambda_p(x,y) + y^{[p]} \quad\text{for all $x,y\in \lieg$},
$$ 
where $\Lambda_p(x,y)$ is a universal expression in terms of iterated Lie brackets
due to Jacobson. Rather than explaining this  tricky matter, 
I refer to \cite{Demazure; Gabriel 1970}, Chapter II, Section 7.2.
Note, however, that in characteristic  $p=2$ things simplify  and we have $(x+y)^{[2]}=x^{[2]} - [x,y] + y^{[2]}$.

The group scheme $\alpha_p$ is entirely determined by
its restricted Lie algebra $\Lie(\alpha_p)$.  This
correspondence works, more generally,    for finite infinitesimal group schemes $G$
of height $\leq 1$.
Let me recall the basic steps in this correspondence.
If $\lieg=\Lie(G)$ is the restricted Lie algebra,
one may interpret the restricted universal enveloping algebra $U^{[p]}(\lieg)$
as the \emph{algebra of distributions} on $G$ at the neutral element $e\in G$, and its $k$-linear dual $U^{[p]}(\lieg)^\vee$
as the \emph{algebra of functions} on $G$. In particular, we have
$$
G=\Spec(U^{[p]}(\lieg)^\vee).
$$
Here the  diagonal  map $\lieg\ra \lieg\oplus \lieg$
induces  the multiplication in $U^{[p]}(\lieg)^\vee$.
Recall that $U^{[p]}(\lieg)$ is the quotient
of the universal enveloping algebra $U(\lieg)$ modulo the ideal generated by
the elements $x^p-x^{[p]}$ with $x\in\lieg$. 
From this discussion it follows that the canonical
map
$$
\Hom(G,H)\lra\Hom(\lieg,\lieh),\quad \varphi\longmapsto\Lie(\varphi)
$$
is bijective for all group schemes $H$, where $\lieh=\Lie(H)$.

Now consider the special case that $H=\Aut_{Y/k}$ for some $k$-scheme $Y$,
such that $\lieh=H^0(Y,\Theta_{Y/k})$. Then the set of $G$-actions on $Y$ is
in correspondence to the set of homomorphisms $\varphi:\lieg\ra H^0(Y,\Theta_{Y/k})$ of restricted Lie algebras.
Indeed, such a homomorphism of restricted Lie algebras induces a homomorphism of $k$-algebras
$U^{[p]}(\lieg)\ra\Diff(\O_Y,\O_Y)$ into the algebra of differential operators, and the adjoint map
$\O_Y\ra\O_Y\otimes_k U^{[p]}(\lieg)^\vee$ then defines the desired action
$ Y\times G\ra Y$.

In the   case  $G=\alpha_p$, the preceding discussion simplifies as follows. To give a homomorphism
of group schemes $\varphi:\alpha_p\ra H$
is nothing but to give a vector $\delta\in\Lie(H)$ with $\delta^{[p]}=0$,
by setting $\delta=\Lie(\varphi)(\partial)$.
As a special case we have $\End(\alpha_p)=k$ and $\Aut(\alpha_p)=k^\times$.
If we change the vector $\delta$ by a nonzero factor $\lambda\in k$,
we only compose the homomorphism $\varphi:\alpha_p\ra H$ with an automorphism of $\alpha_p$, but do not change its image.
Note that the subset
$$
\left\{\delta\in\Lie(H)\mid \delta^{[p]}=0\right\}\subset \Lie(H)
$$
is a cone, but in general not a subgroup. It follows that the space
of nonzero homomorphisms $\varphi:\alpha_p\ra H$ may have several connected components.
For later use we record:

\begin{proposition}
\mylabel{action vector}
Let $Y$ be a $k$-scheme. Then there is a bijection between
the set of $\alpha_p$-actions on $Y$ and the vectors $\delta\in H^0(Y,\Theta_{Y/k})$ with
$\delta^{[p]}=0$.
\end{proposition}

It is worthwile to see  this correspondence explicitly: Write $G=\Spec k[t]/(t^p)$ and 
$\lieg=k\partial$, where $\partial(t^i)=\delta_{i1}$ is given by the Kronecker delta. 
The restricted universal enveloping algebra is the truncated polynomial algebra
$U^{[p]}(\lieg)=k[\partial]/(\partial^p)$.
The monomials $\partial^i$ with $0\leq i<p$ form a basis in $U^{[p]}(\lieg)$. Its dual basis in 
$U^{[p]}(\lieg)^\vee$ are the divided powers $\gamma_i=t^i/i!$.
Now suppose we have  a vector field $\delta\in H^0(Y,\Theta_{Y/k})$ with $\delta^{[p]}=0$.
Then the corresponding $\alpha_p$-action on $Y$ is given by the formula
\begin{equation}
\label{action formula}
\O_Y\lra\O_Y\otimes U^{[p]}(\lieg)^\vee,\quad
f\longmapsto \sum_{i=0}^{p-1}\frac{\delta^i(f)}{i!}\otimes t^i.
\end{equation}
Note that the formula involves differential operators $\delta^i$, $1<i<p$, which are
in general not derivations.
For $p=2$, however, the situation simplifies a lot, for we just have $f\mapsto f\otimes 1+\delta(f)\otimes t$.

Using Formula (\ref{action formula}), it is now easy to understand
invariant or fixed   subschemes:
Let $A\subset Y$ be a closed subscheme with ideal $\shI\subset \O_Y$.
It follows from Formula (\ref{action formula}) that $A$ is \emph{invariant} if and only if
$\delta(\shI)\subset \shI$, or equivalently that the composition $\delta:\shI\ra\O_A$ is zero.
The schematic image $\alpha_p\cdot A\subset Y$ of the   morphism $G\times A\ra Y$
is called the \emph{orbit} of $A$. Its ideal $\shJ\subset\O_Y$ is the intersection of the kernels
for the $k$-linear maps $\delta^i:\shI\ra\O_A$ for $i=1,\ldots,p-1$, again by Formula (\ref{action formula}).
For $p=2$ we simply have $\shJ=\ker(\delta:\shI\ra\O_A)$.

The induced action on an invariant closed subscheme $A\subset Y$ is trivial if and only if
$\delta(\O_Y)\subset\shI$. Equivalently, the composite map $\delta:\O_Y\ra\O_A$ vanishes.
In particular, a closed point $y\in Y$  with maximal ideal $\maxid\subset\O_{Y}$
is a \emph{fixed point} if and only
if $\delta(\O_{Y})\subset\maxid$.
The fixed closed subsets correspond to
the quasicoherent ideals containing the abelian subsheaf $\delta(\O_Y)\subset\O_Y$.

\section{Quotients, singularities, and base change}
\mylabel{quotients}

We keep the situation as in the preceding section.
Let $\delta\in H^0(Y,\Theta_{Y/k})$ be a vector field with 
$\delta^{[p]}=0$, which defines an $\alpha_p$-action on $Y$.
In this section we discuss quotients and their properties.
 According to \cite{Demazure; Gabriel 1970}, Chapter III, Proposition 3.2
quotients by finite infinitesimal group schemes always exist.
For $\alpha_p$ this works as follows:
The underlying topological space for $Z=Y/\alpha_p$ is homeomorphic to $Y$, so that
the quotient map $Y\ra Z$ is a homeomorphism, and the 
structure sheaf $\O_Z$ is the kernel of the derivation $\delta:\O_Y\ra\O_Y$.
If $Y$ is of finite type over $k$, so is $Z$, and the morphism $Y\ra Z$ is finite.
In any case,  this morphism is integral. 
From this it follows that there is a unique maximal fixed closed subset   $Y^{\alpha_p}\subset Y$,
which is called the \emph{fixed scheme}. Its ideal is the quasicoherent ideal that is locally generated
by the abelian subsheaf $\delta(\O_Y)\subset\O_Y$.
If there are no fixed points, then $Y\ra Z$ is an $\alpha_p$-torsor in the fppf-topology,
by \cite{Demazure; Gabriel 1970}, Chapter III, Proposition 3.2.
In particular, $Y\ra Z$ is flat of degree $p$, with  geometric fibers   isomorphic
to the spectrum of $\bar{\k}(z)[t]/(t^p)$.

We now discuss singular and nonsmooth points on the quotient.
Let $y\in Y$ be a point, and $z\in Z$ be its image.
We assume that $Y$ and hence $Z$ are of finite type over $k$. The following
two results are very useful in determining the nonsmooth locus on the quotients $Z$.

\begin{proposition}
\mylabel{isolated fixed}
Suppose that $Y$ is of dimension $\geq 2$ and satisfies Serre's Condition $(S_2)$.
Let $y\in Y$ be an isolated fixed point.
Then $Z$ is not smooth near  $z\in Z$.
\end{proposition}

\proof
The problem is local in $Z$. Seeking a contradiction, we assume that $Z$ is affine and smooth,
and that $y$ is the only fixed point.
Set  $U=Y-\left\{y\right\}$ and $V=Z-\left\{z\right\}$. Then the  quotient morphism $Y\ra Z$
induces an  $\alpha_p$-torsor $U\ra V$. According to Ekedahl's purity result \cite{Ekedahl 1988}, Proposition 1.4,
this $\alpha_p$-torsor extends to an $\alpha_p$-torsor $T\ra Z$.
By assumption, $Y$ satisfies Serre's condition $(S_2)$, and $T$ is obviously Cohen--Macaulay.
It follows that both restriction maps 
$H^0(T,\O_T)\ra H^0(U,\O_T)$ and $H^0(Y,\O_Y)\ra H^0(U,\O_Y)$ are bijective,
which implies $Y= T$. The maps $\O_Y\ra\O_Y\otimes k[t]/(t^p)$ and
$\O_T\ra\O_T\otimes k[t]/(t^p)$ defining the  
$\alpha_p$-actions  are uniquely determined by their restriction $U$. 
It follows that the action on $Y$ is free, contradiction.
\qed

\medskip
Over nonperfect ground fields $k$, one has to distinguish regularity and geometric regularity (= formal smoothness).
I do not know if, in the preceding situation, the quotient $Z$ could be regular.
The next lemma   deals with regularity rather then smoothness.

\begin{proposition}
\mylabel{projective dimension}
Suppose   $y\in Y$ is not a fixed point. Let $\shI\subset\O_Y$ be the ideal of
the orbit $A=\alpha_p\cdot \left\{y\right\}$.
Then $\O_{Z,z}$ is regular if and only if the ideal $\shI \subset\O_{Y,y}$ has
finite projective dimension. In this case, $\shI $ is generated by $n=\dim(\O_{Y,y})$ elements.
\end{proposition}

\proof
The problem is local in $Z$, so we may assume that $Z$ is affine and that
$Y\ra Z$ is an $\alpha_p$-torsor.
Suppose first  that $\O_{Z,z}$ is regular. Then the maximal ideal $\maxid\subset\O_{Z,z}$ has finite projective
dimension, and is generated by $n$ elements.
The exact sequence $0\ra\maxid\ra\O_Z\ra\k(z)\ra 0$
induces by flatness an exact sequence
$$
0\lra \maxid\otimes\O_Y\lra\O_Y\lra\O_A\lra 0
$$
hence $\maxid\otimes\O_Y=\shI$. Using flatness again, we infer that 
$\shI $ has finite projective dimension and is generated by
$n$ elements.

Suppose conversely $\shI $ has finite projective dimension, say $\pd(\shI )=m$.
Choose a resolution
$$
0\lra M\lra F_m\lra \ldots\lra F_0\lra\maxid\lra 0
$$
with $F_0,\ldots, F_m$ free and finitely generated.
Pulling back, we obtain an exact sequence
$$
0\lra M\otimes\O_{Y,y}\lra F_m\otimes\O_{Y,y}\lra \ldots\lra F_0\otimes\O_{Y,y}\lra\shI \lra 0
$$
By Hilbert's Syzygy Theorem, $M\otimes\O_{Y,z}$ is free.
It follows from descent theory that already $M$ must be free (\cite{SGA 1}, Expos\'e VII, Corollary 1.11).
Hence $\pd(\k(z))<\infty$, so the local ring $\O_{Z,z}$ is regular.
\qed

\medskip
This gives a handy criterion in terms on embedding dimensions on $Y$ for singularities on $Z$:

\begin{corollary}
\mylabel{embedding criterion}
Suppose that $y\in Y$ is a rational point that is not fixed, with
$\edim(\O_{Y,y})\geq \dim(\O_{Y,y}) +p$. Then the local ring $\O_{Z,z}$ is not regular.
\end{corollary}

\proof
Let $\maxid\subset\O_{Y,y}$ be the maximal ideal of $y\in Y$, and
$ I\subset\O_{Y,y}$ be the ideal for its orbit.
Then we have a short exact sequence $0\ra I\ra\maxid\ra \maxid/I\ra 0$, where
the quotient $\maxid/I$ is a $k$-vector space of dimension $p-1$.
Tensoring with $k$ gives an exact sequence
$$
\Tor_1(\maxid,k)\lra \Tor_1(k,k)\lra  I\otimes k\lra\maxid\otimes k \lra \maxid/I\lra 0.
$$
Whatever the contribution from the Tor terms on the left, the $k$-vector space
$I\otimes k$ has dimension $\geq\edim(\O_{Y,y})-(p-1)\geq \dim(\O_{Y,y})+1$. By the Nakayama Lemma, it is impossible
to generate $I$ with $n=\dim(\O_{Y,y})$ elements.
According to Proposition \ref{projective dimension}, the local ring $\O_{Z,z}$ must be singular.
\qed

\medskip
Provided that the group scheme action is free, 
the scheme $Y$ satisfies Serre's condition $(S_n)$ if and only if
the quotient $Z$ satisfies   $(S_n)$, by \cite{EGA IVb}, Corollary 6.4.2.
For nonfree action, we at least have the following:

\begin{proposition}
\mylabel{serre condition}
Suppose that $Y$ satisfies Serre's condition $(S_2)$.
Then the quotient $Z=Y/\alpha_p$ satisfies Serre's condition $(S_2)$ as well.
\end{proposition}

\proof
We  may assume that $Z$ is affine.
Let $V\subset Z$ be an open subset
whose complement has codimension $\geq 2$.
We have to check that the restriction map $H^0(Z,\O_Z)\ra H^0(V,\O_Z)$ is bijective.
Injectivity is clear, because $\O_Z\subset\O_Y$ contains no torsion sections.
As to surjectivity, suppose $s_V\in H^0(V,\O_Z)$.
By assumption, it extends to a section $s\in\Gamma(Y,\O_Y)$. The section
$\delta(s)\in H^0(Y,\O_Y)$ vanishes on $V$, hence vanishes everywhere, whence
$s$ is a section for $\O_Z$ extending $s_V$.
\qed

\medskip
The following result on   dualizing sheaves on quotients is useful in
constructing Calabi--Yau quotients; we shall use it in Section \ref{k3}.

\begin{proposition}
\mylabel{dualizing quotient}
Suppose that $Y$ is Gorenstein with $\omega_Y=\O_Y$, that the fixed scheme $Y^{\alpha_p}$ is empty, and
that $k=\Gamma(Y,\O_Y)$. Then $Z$ is Gorenstein with $\omega_Z=\O_Z$.
\end{proposition}

\proof
The morphism $q:Y\ra Z$ is flat with Gorenstein fibers, because the group scheme action is free.
Hence $\omega_{Y/Z}$ is invertible. By assumption, the dualizing sheaf 
$\omega_{Y}=\omega_{Y/Z}\otimes q^*(\omega_Z)$
is invertible. Using descent theory, we infer that $\omega_Z$ is invertible, that is,
$Z$ is Gorenstein.

We next check that $\omega_{Y/Z}$ is trivial. Let $\shI\subset\O_{Y\times_Z Y}$ be the
ideal of the diagonal. The group scheme action is free by assumption, hence
$Y\times_ Z Y=Y\times\alpha_p$. It follows that $\shI\simeq\O_Y$ as $\O_Y$-module.
According to Kunz \cite{Kunz 1986}, page 363 there is a canonical map
$\shI\ra\shHom_{\O_Y}(\shHom_{\O_Z}(\O_Y,\O_Z),\O_Y)$, which is bijective.
The term on the right is the dual of $\omega_{Y/Z}$, which therefore is trivial.

It remains to check that the preimage map $\Pic(Z)\ra\Pic(Y)$ is injective.
For this we may assume that our ground field $k$ is algebraically closed.
Jensen showed in \cite{Jensen 1978}, Section 2 that we have an exact sequence
$$
0\lra D(\alpha_p)\lra \underline{\Pic}_{Z/k}\lra \underline{\Pic}_{Y/k}.
$$
of group-valued functors. Here $D(\alpha_p)$ is the Cartier dual, which in our case is
uncanonically isomorphic to $\alpha_p$, hence contains no point but the origin.
Note that Jensen assumed that $Y$ is proper, but his arguments go through with
the weaker assumption $k=\Gamma(Y,\O_Y)$.
We infer that $\Pic(Z)\ra\Pic(Y)$ is injective,  and hence $\omega_Z=\O_Z$.
\qed

\medskip
We next discuss base-change properties for quotients.
Suppose that $Y$ is endowed with a morphism $h:Y\ra \Spec(R)$ onto some affine $k$-scheme,
and that our derivation $\delta$ lies in $H^0(Y,\Theta_{Y/R})\subset H^0(Y,\Theta_{Y/k})$.
Then the group scheme action $\alpha_p\times Y\ra Y$ is an $R$-morphism,
and the quotient $Z$ is an $R$-scheme. 
The following base-change
property is standard:

\begin{proposition}
\mylabel{base free}
The formation of $Z=Y/\alpha_p$ commutes with flat base change in $R$.
If   the fixed scheme $Y^{\alpha_p}$ is empty, it commutes with arbitrary base change in $R$.
\end{proposition}

\proof
 Consider the exact sequence of quasicoherent  $\O_Z$-modules
\begin{equation}
\label{base change}
0\lra \O_Z\lra\O_Y\stackrel{\delta}{\lra}\O_Y
\end{equation}
on $Z$. Given  a flat $R$-algebra $R'$, we obtain another exact sequence by applying $\otimes_R R'$,
so taking the quotient commutes with flat base change.

Now suppose our action is free.
The inclusion $\O_Z\subset \O_Y$ remains injective after tensoring with $\k(z)$ for
all $z\in Z$, because $1\in\O_{Z,z}$.
Since the $\alpha_p$-action is free, the quasicoherent $\O_Z$-module $\O_Y$ is locally free
of rank $p$. Shrinking $Z$, we may assume that the inclusion admits a splitting, hence
the the exact sequence (\ref{base change}) remains exact after tensoring with arbitrary $R$-algebras $R'$.
\qed

\medskip
For the geometric constructions I have in mind, it is   crucial to
work with group actions with fixed points.  I see no   reason
why  the  base-change property should   hold in general. To discuss this problem, suppose
for simplicity that $Y=\Spec(A)$ is affine, such that
$Z=\Spec(B)$, where $B$ is defined by the exact sequence
\begin{equation}
\label{flat sequence}
0\lra B\lra A \stackrel{\delta}{\lra} A\lra \coker(\delta)\lra 0.
\end{equation}
The issue is as follows:

\begin{proposition}
\mylabel{flatness discussion}
Suppose that $Y$ is $R$-flat. Then the quotient
$Z$ is $R$-flat if and only if the $R$-module $\coker(\delta)$ has flat
dimension $\leq 2$. The formation of $Z=Y/\alpha_p$ commutes with arbitrary
base change in $R$ if and only if $\coker(\delta)$ is $R$-flat.
\end{proposition}

\proof
By assumption, the $R$-module $A$ is flat. 
The exact sequence  (\ref{flat sequence})
gives 
$\flatdim(B)=\flatdim(\coker(\delta))-2$ for flat dimensions, hence the first assertion.
For the second assertion,  let $R'$ be an $R$-algebra, and  consider the induced sequence
$$
0\lra B\otimes R'\lra   A\otimes R' \stackrel{\delta\otimes 1}{\lra} A\otimes R'.
$$
Suppose that $\coker(\delta)$ is flat. Then the exact sequence (\ref{flat sequence})
is a flat resolution, so the preceding sequence remains exact,
which means that the formation of the quotient commutes with base change.
Conversely, suppose the preceding sequence stays exact for all $R'$.
Then $\Tor_1^R(\coker(\delta),R')=0$, whence $\coker(\delta)$ is flat.
\qed

\medskip
Let me record the following consequence:

\begin{corollary}
\mylabel{2-dimensional flat}
Suppose that $R=k[r,s]$ is a polynomial ring in two variables.
Then $Z\ra\Spec(R)$ is flat. 
\end{corollary}

\proof
The ring $R=k[r,s]$ has homological dimension
$\homdim(R)=2$, hence also Tor-dimension $\Tordim(R)=2$
(confer  \cite{Weibel 1994}, Section 4.1).
Hence $\coker(\delta)$ has projective dimension $\leq 2$, so Proposition \ref{flatness discussion} applies.
\qed

\medskip
We shall use the following    base-change-property in Section \ref{selfproduct}, which I formulate
for a rather special situation:
Suppose that $p=2$ and  $A=R[x,y]$ is a polynomial algebra in two variables.
Then we have $\delta=f D_x + gD_y$ for some   $f,g\in R[x,y]$,
where $D_x=\partial/\partial x$ and $D_y=\partial/\partial_y$.

\begin{corollary}
\mylabel{2-dimensional change}
In the preceding situation, suppose that $f,g\in R[x,y]$ are monic polynomials.
Then the formation of $Z=Y/\alpha_2$ commutes with arbitrary base change in $R$.
\end{corollary}

\proof
We check that $\coker(\delta)$ is a free $R$-module.
All terms in the exact sequence (\ref{flat sequence}) are modules
over $A'=R[x^2,y^2]$, because the derivation $\delta$
is $A'$-linear. Moreover, $A$ is a free $A'$-module of rank four,
with basis $1,x,y,xy$. The matrix of $\delta$ with respect to this basis
is 
$$
\begin{pmatrix}
0 & f & g & 0\\
0 & 0 & 0 & g\\
0 & 0 & 0 & f\\
0 & 0 & 0 & 0\\
\end{pmatrix},
$$
hence $\coker(\delta)=  (A'/fA' + gA') \oplus A'^{\oplus 2}/(f,g)A'\oplus A'$.
The first summand is a free $R$-modules, because $f,g$ are monic,
and the last summand is obviously free.
The middle summand sits in an exact sequence
$$
0\lra A'\lra A'^{\oplus 2}/(f,g)A' \stackrel{pr_2}{\lra} A'/gA'\lra 0,
$$
hence is a free $R$-module as well.
\qed

\section{The cuspidal rational curve}
\mylabel{cuspidal}

Let $k$ be a ground field of characteristic $p=2$, and
$C$ be the cuspidal curve of arithmetic genus $p_a(C)=1$.
To be explicit, we choose an indeterminate $u$ and  write
$$
C=\Spec k[u^2,u^3] \;\cup\; \Spec k[u^{-1}].
$$
The goal of this section is to determine the restricted Lie algebra
of $\Aut_{C/k}$ and read off all possible $\alpha_2$-actions on $C$.
This extends results of Bombieri and Mumford (\cite{Bombieri; Mumford 1976}, Proposition 6),
who already studied $\Aut_{C/k}$ in connection with
quasielliptic fibrations in the Enriques classification of surfaces.

I start by describing the   Lie algebra  $\lieg=H^0(C,\Theta_{C/k})$ in   Lie-theoretic terms
without referring to geometry:
Let $\liea$ be the 3-dimensional restricted Lie algebra over $k$ with trivial
Lie brackets $[a,a']=0$ and trivial $p$-map $a^{[p]}=0$. Let $\lieb$ be
the 1-dimensional restricted Lie algebra endowed with a generator $b\in \lieb$
with $b^{[p]}=b$. Any linear endomorphism of $\liea$ is a derivation, because $\liea$ is
commutative. The linear map
$\lieb\ra\liegl(\liea)$, $b\mapsto \id_\liea$
commutes with Lie brackets.
As explained in \cite{Bourbaki LIE 1}, \S 1.7 this homomorphism yields a semidirect product Lie algebra
$\lieg=\liea\rtimes\lieb$,
whose Lie bracket   is  
$$
[a+\lambda b,a'+\lambda'b]=\lambda a'-\lambda' a.
$$
We remark in passing that the ideal $\liea\subset\lieg$ is the \emph{derived ideal} $\lieD\lieg=[\lieg,\lieg]$,
whence $\lieD^2\lieg=0$ and $\lieg$ is solvable. 
On the other hand, $\lieg$ has trivial center, so $\lieg$ it is not nilpotent.

\begin{lemma}
\mylabel{unique restricted}
There is precisely one $p$-map on the semidirect product Lie algebra $\lieg=\liea\rtimes\lieb$
such that the canonical inclusions $\liea,\lieb\subset\lieg$ are restricted inclusions.
\end{lemma}

\proof
If it exists, the $p$-map for $\lieg$ is uniquely determined:
For $p=2$ the axioms for $p$-maps give
$(x+y)^{[2]}=x^{[2]}+y^{[2]} - [x,y]$. Therefore the $p$-map   
in our semidirect product must be
\begin{equation}
\label{power map}
(a+\lambda b)^{[2]}=a^{[2]} + (\lambda b)^{[2]} + [a,\lambda b] =
\lambda (a+\lambda b).
\end{equation}
It remains to check that this   meets the remaining two axioms for $p$-maps.
Indeed, given any scalar $\mu\in k$, we obviously have
$$
(\mu(a+\lambda b))^{[2]}=\mu\lambda(\mu(a+\lambda b)) = \mu^2(a+\lambda b)^{[2]}.
$$
Moreover, we easily compute
$$
[(a+\lambda b)^{[2]},a'+\lambda'b] = [\lambda a+\lambda^2b,a'+\lambda'b] = \lambda^2a'+\lambda\lambda'a,
$$
which equals
$$
[a+\lambda b,[a+\lambda b,a'+\lambda'b]]=[a+\lambda b,\lambda a'+\lambda' a] = \lambda(\lambda a'+\lambda'a).
$$
Hence Formula (\ref{power map}) indeed defines a $p$-map.
\qed

\medskip
The book of Strade and Farnsteiner \cite{Strade; Farnsteiner 1988} contains
more information on the existence of $p$-maps in Lie algebras.
We now come back to our cuspidal curve $C$ of arithmetic genus one:

\begin{proposition}
The restricted Lie algebras $H^0(C,\Theta_{C/k})$ and $\lieg=\liea\rtimes\lieb$
are isomorphic.
\end{proposition}

\proof
Consider the two affine open subsets $U=\Spec k[u^2,u^3]$  
and $V=\Spec k[u^{-1}]$.
The only relation between the generators $u^2,u^3$ on $U$ is the obvious one, namely $(u^2)^3=(u^3)^2$.
Whence $\Omega^1_{U/k}$ is generated by the differentials $d(u^2),d(u^3)$ modulo the
relation $u^4d(u^2)=0$. Since $\O_C$ is torsion free, the dual $\Theta_{U/k}$ is a free
$\O_U$-module of rank one, generated by the form $d(u^3)\mapsto 1$. 
On the overlap $U\cap V$, this form becomes $u^{-2}D_u$, where
$D_u=\partial/\partial u$ is taking derivative with respect to $u$. We shall use the same  symbol
$u^{-2}D_u$ to denote
this form on $U$, although the two individual factors $u^{-2}$ and $D_u$ do not make sense on $U$.

The affine open subset $V$ is smooth.
Here $\Omega^1_{V/k}$ is free of rank one, generated by the differential
$d(u^{-1})$. The dual $\Theta_{V/k}$ is free of rank one as well, generated by
the form $d(u^{-1})\mapsto 1$.
On the overlap $U\cap V$ we have $d(u^{-1})=u^{-2}du$, so our form becomes
$u^2D_u$. From this we infer that the Lie algebra $\lieg'=H^0(C,\Theta_{C/k})$ is a $k$-vector space
of rank 4, with basis
$$
u^{-2}D_u,\quad D_u,\quad uD_u,\quad u^2D_u.
$$
Let $\liea\subset\lieg'$ be the subspace generated by the derivations
with even coefficients, that is, $u^{-2}D_u,D_u,u^2D_u$.
Then one easily computes that both   Lie bracket and   $p$-map
vanish on $\liea$. Moreover, we have $[uD_u,u^{2i}D_u]=u^{2i}D_u$ for all integers $i$.
Whence $\liea\subset\lieg'$ is the derived Lie algebra for $\lieg'$.
The resulting extension of   Lie algebras
$$
0\lra\liea\lra\lieg'\lra\lieb\lra 0
$$
splits, because $\lieb$ is 1-dimensional. The element $uD_u\in\lieg'$ defines a splitting, 
and we have $(uD_u)^{[2]}=uD_u$.
We now   regard $\lieg'$ as a semidirect product.
Any such semidirect product is given by a homomorphism $\lieb\ra\liegl(\liea)$.
In our case, this map sends the generator $uD_u$ to the identity $\id_{\liea}$.
It follows that $\lieg'=H^0(C,\Theta_{C/k})$ and $\lieg=\liea\rtimes\lieb$
are isomorphic as Lie algebras. By Lemma \ref{unique restricted}, they also have
the same $p$-maps.
\qed

\medskip
Our algebraic computations translate into the following geometric statement: 

\begin{corollary}
\mylabel{parameter space}
The set of $\alpha_2$-operations on the cuspidal curve $C$ of arithmetic genus one is parameterized by
the affine space $\AA^3$.
\end{corollary}

\proof
According to Proposition \ref{action vector}, the set in question is the set of all $a+\lambda b\in \lieg$
with $(a+\lambda b)^{[2]}=0$. By formula (\ref{power map}), each vector $a$ from the
derived ideal $\liea\subset\lieg$  has this property.
On the other hand, each vector $a+\lambda b\in\lieg$ with $\lambda\neq 0$ has
$(a+\lambda b)^{[2]}=\lambda(a+\lambda b)\neq 0$.
\qed

\section{The selfproduct for the cuspidal rational curve}
\mylabel{selfproduct}

 Let $C$ be 
the cuspidal curve of arithmetic genus one over a ground field $k$ of characteristic $p=2$,
as in the preceding  Section.
Throughout we consider the selfproduct 
$Y=C\times C$, which is a nonnormal integral proper surface
with unibranch singularities and normalization $\PP^1\times\PP^1$.
We shall use coordinates $u^2,u^3,u^{-1}$ on the first factor of $C\times C$, and coordinates
$v^2,v^3,v^{-1}$ on the second factor:
$$
C=\Spec k[u^2,u^3] \cup \Spec k[u^{-1}] \quadand C=\Spec k[v^2,v^3] \cup \Spec k[v^{-1}].
$$
I want  to define $\alpha_2$-actions on $Y=C\times C$ depending on two parameters $r,s$.
In the following, we have to allow parameters from various parameter rings
(for example fields, dual numbers, polynomial rings and so forth).
To unify notation, let us fix a $k$-scheme $S$ and two global sections  $r,s\in H^0(S,\O_S)$,
and consider the scheme
$$
Y\times S=(C\times C)\times S,
$$
viewed as a proper flat family of surfaces over $S$.
The derivation
$$
\delta=(u^{-2}+r)D_u + (v^{-2} + s)D_v
$$
defines a global vector field $\delta\in   H^0(Y\times S,\Theta_{Y\times S/S})$. 
A straightforward computation shows that $\delta\circ \delta=0$. Hence $\delta$
defines  
an $\alpha_{2}$-action on $Y\times S$ over $S$.
Let us first examine the fixed scheme   for this action.

\begin{proposition}
\mylabel{fixed scheme}
The fixed scheme $(Y\times S)^{\alpha_2}\subset Y\times S$ is contained in the open subset
$\Spec \O_S[u^{-1},v^{-1}]\subset Y\times S$, and its ideal   is generated by the two elements
$u^{-2}(u^{-2}+r)$ and $v^{-2}(v^{-2}+s)$.
\end{proposition}

\proof
The ideal $\shI\subset\O_{Y\times S}$ generated by the image of the derivation $\delta$
defines the fixed scheme in question.
Since $\delta(u^3)=1+ru^2$, we have $\shI_y=\O_{Y,y}$ for all points $y\in Y\times S$ outside $\Spec \O_S[u^{-1},v^{-1}]$.
Over $\Spec \O_S[u^{-1},v^{-1}]$,  we compute 
\begin{equation}
\label{inverse derivation}
\delta=u^{-2}(u^{-2}+r)D_{u^{-1}} + v^{-2}(v^{-2}+s)D_{v^{-1}},
\end{equation}
and the result follows.
\qed

\medskip
We see that the geometric fibers of $Y\times S\ra S$ contain precisely four fixed points
over the open subset $D(rs)\subset S$. These fixed points come together in pairs over
the closed subsets $V(r)$ and $V(s)$.
On the intersection $V(r,s)=V(r)\cap V(s)$, precisely one fixed point remains.

Next, we  examine the quotient scheme $\foZ=Y\times S/\alpha_2$ with respect to the $\alpha_2$-action on 
the product family $Y\times S$. Here I use the fracture letter $\foZ$, to emphasize that
the resulting morphism $\foZ\ra S$ usually is not a product family
(this has nothing   to do with formal schemes).

\begin{proposition}
\mylabel{flatness change}
The resulting morphism $\foZ\ra S$ is flat and commutes with arbitrary base change in $S$.
\end{proposition}

\proof
We first check the base-change property.
By Proposition \ref{base free}, the base-change property  holds outside the fixed points.
In light of Corollary \ref{2-dimensional change} and Formula (\ref{inverse derivation}), the
base-change property holds near the fixed points as well.
To check flatness, it therefore suffice to treat the universal situation $R=k[s,t]$,
and then flatness holds by Corollary \ref{2-dimensional flat}.
\qed

\medskip
Throughout the paper, we shall frequently pass back and forth
between the family $\foZ\ra S$ and its fibers $\foZ_\sigma$, $\sigma\in S$.
We just saw that the fiber $\foZ_\sigma$ is also the quotient of $Y_\sigma$
by the induced $\alpha_2$-action. To simplify notation, we usually write $Z=\foZ_\sigma$
to denote fibers; making base change, we then usually assume that $S$
is the spectrum of our ground field $k$ and write $Z=\foZ$.

\begin{proposition}
\mylabel{family normal}
The fibers $Z=\foZ_\sigma$ of the flat family $\foZ\ra S$ are normal.
\end{proposition}

\proof
We may assume that $S=\Spec(k)$ and that $k$ is algebraically closed.
By Proposition \ref{serre condition}, the surface $Z$ is Cohen--Macaulay.
It remains to check that it is regular in codimension one.
For this  we may ignore the fixed points, which are isolated by Proposition \ref{fixed scheme}.
Let $U\subset Y$ be the regular locus minus the fixed locus;
by Proposition \ref{projective dimension}, the quotient $Z$ is regular on
the image of $U$.
To finish the argument, let $y\in Y$ be a singular point defined by 
the maximal ideal $\maxid=(u^2,u^3,v^{-1}+\lambda)$ for some $\lambda\in k$
with $\lambda\neq 0,s$. Then the ideal
$$
\shI=(u^2, (1+ru^2)(v^{-1}+\lambda)+(v^{-4}+sv^{-2})u^3)
$$
is $\delta$-invariant, whence the spectrum of $\O_Y/\shI=k[u^3]/(u^6)$ must be the orbit of $y$.
Clearly, the ideal $\shI$ has finite projective dimension, whence
  $Z$ is regular at the image of $y$.
Summing up, $Z$ is regular in codimension one.
\qed

\medskip
In light of Proposition \ref{isolated fixed}, the four   fixed points on $Y$
contribute to the singular points $z\in Z$.
It turns out that these singularities are \emph{geometrically} isomorphic,
because the automorphism group scheme $\Aut_{Z/k}$ acts transitively on them. More precisely:

\begin{proposition}
\mylabel{automorphism group}
Assume that $S=\Spec(k)$, and that $r,s\in k$ are squares.
Let $y_1,y_2\in Y$ be   fixed points, and
$z_1,z_2\in Z$ be the corresponding singular points.
Then there are automorphisms $\phi:Y\ra Y$ with $\phi(y_1)=y_2$
and $\psi:Z\ra Z$ with $\psi(z_1)=z_2$ such that the diagram
$$
\begin{CD}
Y @>\phi>> Y\\
@VVV @VVV\\
Z @>>\psi> Z
\end{CD}
$$
commutes.
\end{proposition}

\proof
In case $r=s=0$, there is precisely one fixed point $y\in Y$, and we have
nothing to prove.
Suppose $r\neq 0$.
Let $\phi:Y\ra Y$ be the involution defined by $u^{-1}\mapsto u^{-1} + \sqrt{r}$.
It is easy to see that  the corresponding $\ZZ/2\ZZ$-action commutes with
the $\alpha_2$-action on $Y$.
It therefore descends  to an involution
$\psi:Z\ra Z$.
Now suppose that $s\neq 0$.
Then we have    similar involutions on $Y$ and  $Z$ induced by  $v^{-1}\mapsto v^{-1} + \sqrt{s}$.
Obviously, the subgroup generated by these  involutions acts transitively on the
set of singularities on $Z$ coming from fixed points on $Y$.
\qed

\section{The singularity coming from the quadruple point}
\mylabel{singularity}

We keep the notation from the preceding section,
but we assume for simplicity that our parameter space $S=\Spec(R)$ is affine. Recall that 
$Y=(C\times C)\otimes R$,
and $\foZ=Y\otimes R/\alpha_2$ is the quotient by the $\alpha_2$-action defined by the vector field
$$
\delta=(u^{-2}+r)D_u + (v^{-2}+s)D_v 
$$
with parameters $r,s\in R$.
Set $A=R[u^2,u^3,v^2,v^3]$, and let $B\subset A$ be the kernel of the derivation
$\delta:A\ra A$.
Then $\Spec(A)\subset Y\otimes R$ is an affine open neighborhood containing fiberwise
the point of embedding dimension four, which is defined by
the ideal $(u^2,u^3,v^2,v^3)\subset A$, and $\Spec(B)\subset \foZ$
is an affine open neighborhood containing fiberwise the corresponding singularities.

\begin{proposition}
\mylabel{invariants embedding}
We have 
$B=R[u^2,v^2,(1+s v^2)u^3 + (1+r u^2)v^3]$
as $R$-subalgebras inside $A$. 
\end{proposition}

\proof
In light of Proposition \ref{flatness change}, we may assume
that $R=k[r,s]$ is a polynomial algebra in two variables.
Consider the $k$-subalgebra $A'=R[u^2,v^2]$ inside $A$.
Obviously, we have $A'\subset B$. Moreover, $A$ is a free $A'$-module of rank four,
with basis $1,u^3,v^3,u^3v^3$. Given $\alpha,\beta,\gamma\in A'$, we compute
$$
\delta(\alpha u^3+\beta v^3 +\gamma u^3v^3)=
\alpha(1+ r u^2) + \beta(1 + s v^2) +\gamma(1+ru^2)v^3+\gamma(1+sv^3)u^3.
$$
Suppose this expression vanishes. Comparing coefficients, we see $\gamma(1+sv^3)=0$
and hence $\gamma=0$.
Using the factoriality of $A'$, we infer that $\alpha$ and $\beta$
are multiples of $1+s v^2$ and $1+r u^2$, respectively.
Hence $A=k[u^2,v^2,(1+s v^2)u^3 + (1+r u^2)v^3]$.
\qed

\medskip
The invariant subring $B\subset A$ also admits a simple description as 
a quotient ring:

\begin{proposition}
\mylabel{quotient embedding}
We have $B=R[a,b,c]/(c^2+ a^3+b^3 + s^2a^3 b^2+ r^2a^2b^3)$.
\end{proposition}

\proof
It suffices to treat the universal situation $R=k[r,s]$. In particular,
we may assume that $R$ is noetherian.
Set $B'=R[a,b,c]/(c^2+ a^3+b^3 + s^2a^3 b^2+ r^2a^2b^3)$.
The mapping
$$
a\longmapsto u^2,\quadand b \longmapsto v^2,\quadand c \longmapsto (1+s v^2)u^3 + (1+r u^2)v^3
$$
clearly induces a surjective homomorphism $B'\ra B$.
To see that it is also injective, it suffices to treat the case $R=k$,
by the Nakayama Lemma.
Being a complete intersection, the ring $B'$ is Cohen--Macaulay.
A straightforward computation of $\Omega^1_{B'/k}$ reveals that
the singular locus  on $\Spec(B')$ is 0-dimensional.
It follows that the 2-dimensional ring $B'$ is integral.
Being a surjection between integral 2-dimensional rings $B'\ra B$, the map in question is bijective.  
\qed

\medskip
Let me now recall some facts from singularity theory.
Suppose $(\O_z,\maxid,k)$ is a   local noetherian ring that is normal, 2-dimensional, with 
resolution of singularities 
$X\ra\Spec(\O_z)$.
One says that $\O_z$ is a \emph{rational double point} if $H^1(X,\O_X)=0$ and $\mult(\O_z)=2$.
Each rational singularity comes along with an  \emph{irreducible root systems},
which  are classified by Dynkin diagrams. This goes   as follows:
Let $E\subset S$ be the reduced exceptional divisor, $E_i\subset E$ be its integral components,
and $V$ the real vector space generated by the $E_i$  endowed with the scalar product
$\Phi(E_i,E_j)=-(E_i\cdot E_j)$. Then the vectors $E_i\in V$ form a root basis for
an irreducible root system.
Note that in characteristic $p=2,3,5$ there are nonisomorphic rational double points
with the same root system, as Artin observed in \cite{Artin 1977}.

When the residue field $k$ is algebraically closed, only   Dynkin diagrams of type $ADE$ appear,
which are precisely Dynkin diagrams with equal root lengths.
Such Dynkin diagrams are also called \emph{simply laced}.
When $k$ is not algebraically closed, Dynkin diagrams with unequal root lengths show up as well.
This means that some field extensions $k\subset H^0(E_i,\O_{E_i})$ are nontrivial.
The degree of the field extensions are determined by  the
Cartan numbers.

Now suppose that the parameter ring $R$ is our ground field $k$, 
and let $z\in Z$ be the singularity defined by the maximal ideal $(a,b,c)\subset B$,
which corresponds to the  closed point $y\in Y$  defined by
the maximal ideal $(u^2,u^3,v^2,v^3)\subset A$, which is the unique point of   embedding dimension four.

\begin{proposition}
\mylabel{d4 b3}
The local ring  $\O_{Z,z}$ is a rational double point.
It is of type $D_4$ if the ground field $k$ contains a third root of unity, 
and of type $B_3$ otherwise.
\end{proposition}

\vspace{1em}
\centerline{\includegraphics{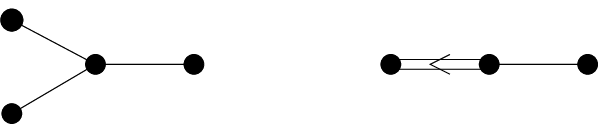}}
\vspace{.8em}
\centerline{ Figure \stepcounter{figure}\arabic{figure}: The Dynkin diagrams $D_4$ and $B_3$.}

\proof
If the ground field $k$ is algebraically closed, it suffices to check for the defining equation
$c^2= a^3+b^3 + s^2a^3 b^2+ r^2a^2b^3$
in one of the available lists of normal forms for rational double points in characteristic two
(see  Artin \cite{Artin 1977} or Greuel and Kr\"oning \cite{Greuel; Kroning 1990}).
To recheck this, and to cover arbitrary ground fields as well, let us do the blowing up $Z'\ra\Spec(B)$
of the ideal $(a,b,c)\subset B$. The scheme $Z'$ is covered by two affine open subsets
$Z'=D_+(a)\cup D_+(b)$. The first one is the spectrum of $k[a,b',c']$, where we set $b'=b/a$ and $c'=c/a$,
modulo  the relation
$$
c'^2=a(1+b'^3) + a^3(s^2+r^2b')b'^2.
$$
A computation with differentials reveals that the reduced locus
of nonsmoothness on $D_+(a)$ lying on the exceptional divisor has ideal $I=(a,1-b'^3)$. 

Now suppose that the ground field $k$ contains a third root of unity $\zeta\in k$.
Then $1-b'^3=(1-b')(1-\zeta b')(1-\zeta^2b')$, hence $D_+(a)$ contains three
singularities, each with residue field $k$, and a straightforward computation
shows that these are rational double points of type $A_1$.
By symmetry, the other open subset $D_+(b)$ does not contribute
any further singularities.
From this we infer that $\O_{Z,z}$ is a rational double point of type $D_4$.

Finally, suppose that $k$ does not contain a third root of unity.
Then we have an irreducible decomposition $1-b'^3=(1-b')(1+b'+b'^2)$,
hence $D_+(a)$ contains two singularities. Again we easily see that
each one is of type $A_1$, but one has residue field $k$, and the
other has residue field isomorphic to the splitting field of $1+b'+b'^2$.
If follows that $\O_{Z,z}$ is a rational double point of type $B_3$.
\qed

\medskip
Finally, suppose that $R$ is a $k$-algebra with residue
field $k$, and set $B_0=B\otimes_R k$ and $R_0=k$.
We now view the $R$-algebra $B$ as a deformation of the $R_0$-algebra $B_0$.
To understand this
deformation near the singularity $z\in \Spec(B_0)$, we have to pass the  
completion
$\hat{B}=R[[a,b,c]]/(c^2+a^3+b^3+s^2a^3b^2+r^2a^2b^3)$.
It turns out that this deformation is trivial:

\begin{proposition}
\mylabel{trivial deformation}
The deformation $\hat{B}$ is isomorphic to the trivial deformation
$\hat{B}_0\hat{\otimes}_{R_0} R$.
\end{proposition}

\proof
Set $f=c^2+a^3+b^3$, and suppose we have a power series $g\in(a^3b^2,a^2b^3)$.
We will show that there is an automorphism $\varphi$ of $R[[a,b,c]]$ with $\varphi(f+g)=f$,
which clearly implies the assertion.

Write $g=\sum\lambda_{ij}a^ib^j$ with $\lambda_{ij}\in R$, and $i,j\geq 2$, and $i+j\geq 5$.
Consider the nonzero monomials in $g$ with minimal total degree $i+j$, and
pick among them the monomial $\lambda_{mn}a^mb^n$ with minimal $a$-degree $m$.
Suppose for the moment $m\geq 3$.
Then $a\mapsto a+\lambda_{mn}a^{m-2}b^n$ defines an automorphism $\varphi_{mn}$ of $R[[a,b,c]]$.
We have 
$$
\varphi(a^3)=a^3+\lambda_{mn}a^mb^n + \lambda^2_{mn}a^{2m-3}b^{2n} + \lambda^3_{mn}a^{3m-6}b^{3n},
$$
and the last two summands have total degree $>m+n$.
A  similar computation of $\varphi(a^k)$, $k\geq 4$ gives that
$g'=\varphi_{mn}(f+g)-f$ lies in $(a^3b^2,a^2b^3)$, and that $g'$ has
only monomials of total degree $>m+n$, or of total degree $m+n$ and $a$-degree $>m$.

In the case $m=2$, we have $n\geq 3$, and we may apply the preceding arguments
with $b$ instead of $a$. Proceeding by induction, we see that we obtain the
desired automorphism $\varphi$ in the form $a\mapsto a +\sum\mu_{ij}a^{i-2}b^j$, $b\mapsto b+\sum\eta_{ij}a^ib^{j-2}$
for certain inductively determined coefficients $\mu_{ij},\eta_{ij}\in R$.
\qed

\begin{remark}
According to \cite{Artin 1977}, there are two isomorphism classes of rational double points
of type $D_4$, which are called $D_4^0$ and $D_4^1$. The upper index has to do with the
versal deformation: The smaller the upper index, the larger the dimension of the versal deformation.
The preceding arguments show that our singularity is of type $D_4^0$.
\end{remark}

\section{Singularities coming from  fixed points}
\mylabel{singularities}

Our next goal is to analyse the singularities on  $\foZ=Y\otimes R/\alpha_2$ coming from
the   fixed points on $Y\otimes R$ for the $\alpha_2$-action.
In this section, we set  $A=R[u^{-1},v^{-1}]$ and define $B\subset A$ to be the kernel
of the derivation $\delta:A\ra A$. Then 
$\Spec(A)\subset Y\otimes R$ is an affine open neighborhood for the fixed points, and
$\Spec(B)\subset \foZ$ is an affine
open neighborhood for    singularities coming from fixed points.

\begin{proposition}
\mylabel{invariants fixed}
We have 
$
B=R[u^{-2},v^{-2}, u^{-1}(v^{-4}+s v^{-2})+v^{-1}(u^{-4}+ru^{-2})]
$
as  subalgebras inside $A$. Fiberwise over $R$, the $B$-module $A$ is reflexive of rank two,
but not locally free.
\end{proposition}

\proof
The arguments for the first statement are as in the proof for Proposition \ref{invariants embedding},
and   left to the reader. For the second statement, we may assume that
$R=k$. 
According to \cite{Serre 1965}, Section IV, Proposition 11, the $B$-module $A$ is Cohen--Macaulay.
Hence  it must be reflexive.
Consider the maximal ideal $$\maxid=(u^{-2},v^{-2},u^{-1}(v^{-4}+s v^{-2})+v^{-1}(u^{-4}+ru^{-2}))$$
inside $B$. We have
$
A/A\maxid = k[u^{-1},v^{-1}]/(u^{-2},v^{-2}),
$
which has length four instead of two. It follows that $A$ is not locally free as $B$-module.
\qed

\medskip
Setting $a=u^{-2}$ and $b=v^{-2}$ and
$c=u^{-1}(v^{-4}+s v^{-2})+v^{-1}(u^{-4}+ru^{-2})$,
we obtain as in Proposition \ref{quotient embedding} the following description of $B$
as a complete intersection:

\begin{proposition}
\mylabel{quotient fixed}
We have $B=R[a,b,c]/(c^2+a(b^4+s^2b^2) + b(a^4+r^2a^2))$.
\end{proposition}

In what follows, we assume $R=k$ and write $Z=\foZ$.
Let $z\in Z$ be the rational point defined by the ideal $\maxid=(a,b,c)$,
which is the image of a fixed point on $Y$.
We want to understand the singularity $\O_{Z,z}$. 
Note that   if the ground field $k$ is perfect, the other  
singularities on $Z$ corresponding to the fixed points on $Y$
are isomorphic to $z\in Z$, by Proposition \ref{automorphism group}.
For our purposes it therefore suffices to understand $\O_{Z,z}$.
We first treat the generic case:

\begin{proposition}
\mylabel{d4}
Suppose  the parameters $r,s\in k$ are both nonzero.
Then the local ring $\O_{Z,z}$ is a  rational double point of type $D_4$.
\end{proposition}

To see this, one has to make an explicit blowing up $X\ra\Spec(B)$ of
the ideal $(a,b,c)$ as in the proof for Proposition \ref{d4 b3}.
We leave this for the reader and immediately turn to the more interesting
case where one of the parameters $r,s$ degenerates:

\begin{proposition}
\mylabel{d8}
Suppose precisely one of the parameters $r,s\in k$ vanishes.
Then the local ring $\O_{Z,z}$ is a rational double point of type $D_8$.
\end{proposition}

\vspace{1em}
\centerline{\includegraphics{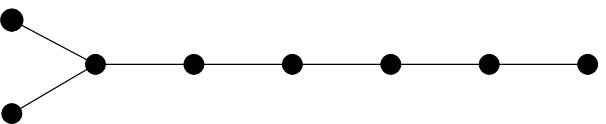}} 
\vspace{.8em}
\centerline{Figure \stepcounter{figure}\arabic{figure}: The Dynkin diagram $D_8$.}
\vspace{1em}

\proof
We may assume $s=0$. Then the defining equation is
$c^2= ab^4 + a^4b + r^2a^2b$
with  $r\neq 0$. We view the corresponding affine scheme
as a double covering of the affine plane $\Spec k[a,b]$.
To produce a resolution for $\O_{Z,z}$ we compute
a blowing-up  $X\ra\Spec k[[a,b]]$
so that the reduced transform of the equation $ab^4 + a^4b + r^2a^2b$
defines a divisor with normal crossings. It turns out that, if  $t$ denotes a local
equation for this normal crossing divisor, the partial derivatives of $t$ vanish precisely
at the singularities of the normal crossing divisor (and luckily nowhere outside the divisor).
Hence the equation $c^2=t$ has at most
rational double points of type $A_1$ (even in characteristic $p=2$), and one immediately reads off the
configuration of the exceptional divisor for a resolution of $\Spec(\O_{Z,z})$.
This is indeed the $D_8$-configuration, and  $X\ra\Spec k[[a,b]]$
is a sequence of three blowing ups.

I do not want to reproduce these computations in detail; but let me explain the first blowing up of $k[a,b]$:
On the $k[a,\frac{b}{a}]$ chart, the transform of $ab^4 + a^4b + r^2a^2b$
is $a^5 \frac{b}{a}+a^5(\frac{b}{a})^4+r^2a^3\frac{b}{a}$.
Removing the square factor $a^2$ we obtain $a\frac{b}{a}(a^2+a^2\frac{b}{a}+r^2)$,
and this is already normal crossing along the exceptional curve $a=0$.

On the $k[\frac{a}{b},b]$ chart, the transform divided by square factors is $(\frac{a}{b})^4b^3+\frac{a}{b}b^3+r^2(\frac{a}{b})^2b$.
There is a unique singularity, which is located at $\frac{a}{b}=b=0$, but this singularity
is not yet normal crossing. Here we have to repeat the process, which I leave to the reader.
\qed

\medskip
It remains to treat the case of totally degenerate parameters $r=s=0$.
Here we do not get   rational singularities. Instead, we get an elliptic singularity.
Let me first recall the  terminology. Suppose that $\O_z$ is
a   local ring that is 2-dimensional and normal, with resolution of singularities $X\ra\Spec(\O_z)$.
The \emph{arithmetic genus} $p_a(\O_z)$ is
defined as the maximal arithmetic genus $p_a(D)=1-\chi(\O_D)$,
where $D\subset X$ ranges over all divisors supported by the
exceptional locus. The rational singularities are precisely those
with $p_a(\O_z)=0$, according to \cite{Artin 1962}, Theorem 1.7. Singularities with $p_a(\O_z)=1$ are called \emph{elliptic}.
Wagreich \cite{Wagreich 1970} obtained a list of elliptic singularities.
For us, the following elliptic singularities are  important:

\vspace{1em}
\centerline{\includegraphics{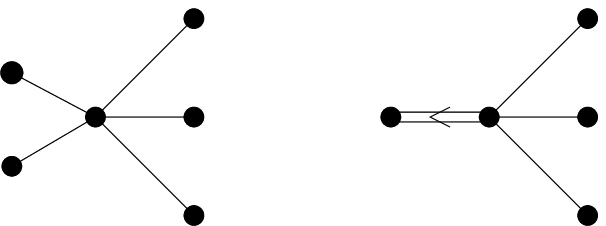}} 
\vspace{.8em}
\centerline{Figure \stepcounter{figure}\arabic{figure}: Wagreich's elliptic singularity of type $19_0$, and a twisted form.}
\vspace{1em}

In the elliptic singularity of type $19_0$, the 
central exceptional curve $E_0$ has selfintersection number $E_0^2=-3$, whereas the
others have $E_i^2=-2$. In the twisted form, the exceptional curve on the left is
$\PP^1_{k'}$ for some   quadratic field extension $k\subset k'$.

\begin{proposition}
\mylabel{elliptic singularity}
Suppose $r=s=0$. Then  $\O_{Z,z}$ is an elliptic singularity.
If the ground field $k$ contains a third root of unity,
it is the elliptic singularity of type $19_0$ from Wagreich's list.
Otherwise, it is the twisted form in Figure 3.
\end{proposition}

\proof
We have $B=k[a,b,c]$ modulo the relation 
$c^2= ab(a-b)(a^2+ab+b^2)$.
If the ground field $k$ contains a third root of unity, the factor $a^2+ab+b^2$
splits into linear factors. According to \cite{Wagreich 1970}, Corollary on page 449,
the singularity is then elliptic of type $19_0$. 
If $k$ does not contain a third root of unity, a similar analysis as in the proof for Proposition \ref{quotient embedding}
shows that the singularity is the twisted form in Figure 3.
\qed

\begin{remark}
Wagreich's elliptic singularity of type $19_0$ also showed up
in Katsura's analysis of the classical Kummer construction in characteristic two \cite{Katsura 1978}.
I do not know whether there is a structural reason for this.
\end{remark}

\section{Quasielliptic fibrations}
\mylabel{quasielliptic}

We keep the notation from  the preceding section and work over a ground field $R=k$, such that
$Z$ is a proper normal surface, defined as the quotient
of $Y=C\times C$ by an $\alpha_2$-action depending on two
parameters $r,s\in k$.
The goal of this section is to study fibrations   on $Z$.

First note that the composite morphism from the normalization
$$
\PP^1\times\PP^1 =\tilde{Y}\lra Y\lra Z
$$
is a finite universal homeomorphism of degree two.
It follows that the induced mapping $\Pic(Z)\ra\Pic(\PP^1\times\PP^1)$ becomes bijective
after tensoring with $\ZZ[1/2]$. Consequently, the proper normal surface $Z$  has Picard number
$\rho(Z)=2$. Moreover, the two projections $Y\ra C$ induce fibrations
$Z\ra\PP^1$, because $C/\alpha_2=\PP^1$.
Throughout, we shall use   the   projection $\pr_2:Y=C\times C\ra C$ onto the second factor,
and denote by $g:Z\ra\PP^1$ the induced fibration, such that the diagram
\begin{equation}
\label{total square}
\begin{CD}
Y @>>> Z\\
@V\pr_2 VV @VVgV\\
C @>>h> \PP^1
\end{CD}
\end{equation}
commutes. In accordance with our   notation $C=\Spec k[v^2,v^3] \cup \Spec k[v^{-1}]$, we write 
the projective line as $\PP^1=\Proj (k[v^2])$. This should cause no confusion.

\begin{proposition}
\mylabel{quasielliptic fibration}
The canonical map $\O_{\PP^1}\ra g_*(\O_Z)$ is bijective, and the
fibration $g:Z\ra \PP^1$ is quasielliptic.
\end{proposition}

\proof
Suppose that the inclusion $\O_{\PP^1}\subset g_*(\O_Z)$ is not bijective.
Using the inclusion $\O_Z\subset\O_Y$   we infer that $g$ factors 
birationally over $C$. It also must factor over the normalization of $\PP^1=\tilde{C}$,
because $Z$ is normal.
In turn, the projection $Y=C\times C\ra C$  
factors over $\tilde{C}$, which is absurd.
Whence $\O_{\PP^1}=g_*(\O_Z)$.

Let $\eta\in\PP^1$ be the generic point, and   $K=\k(\eta)=k(v^2)$ be the corresponding function field.
The generic fiber $Z_\eta=Z_K$ is a curve over the function field $K$. Since $Z$ is normal, $Z_K$ is a regular curve.
Saying that the fibration $f:Z\ra \PP^1$ is \emph{quasielliptic} means that $Z_K$ is a twisted
form of the   cuspidal curve  $C_K$ of arithmetic genus one.
To check this, let $L=k(v)$ be the function field of $C$, and consider the commutative diagram
\begin{equation}
\label{generic square}
\begin{CD}
Y_L @>>> Z_K\\
@VVV @VVV\\
\Spec(L) @>>> \Spec(K).
\end{CD}
\end{equation}
The morphisms are equivariant with respect to the induced $\alpha_2$-actions on $Y_L$ and $\Spec(L)$,
and the horizontal arrows are $\alpha_2$-torsors. More precisely, 
the $\alpha_2$-action on $Y_L=C_K\times\Spec(L)$ is the diagonal action, coming from
actions  on    $C_K$ and $\Spec(L)$.
Passing to the algebraic closure $K\subset\bar{K}$, the induced torsor
$\Spec(L\otimes_K\bar{K})\ra\Spec(\bar{K})$ becomes trivial.
This implies that the quotient of $Y_L\otimes_K \bar{K}$ by the action
is isomorphic to $C_{\bar{K}}$. Summing up,
$Z_{\bar{K}}$ and $C_{\bar{K}}$ are isomorphic.
\qed

\medskip
In what follows, $K\subset L$ denote the function fields for $C\ra\PP^1$, respectively.
Note that, with the notation from the preceding proof, the composite morphism
$$
\PP^1_L\lra Y_L\lra Z_K
$$
is a finite universal homeomorphism.
Our first task is to determine the set of $K$-rational points $Z_K(K)$.
This depends in a strange way on the parameters:

\begin{proposition}
\mylabel{sections}
\begin{enumerate}
\item
Suppose that $r,s\neq 0$, $\frac{s}{r}\in\FF_4$, $\sqrt{r}\in K$.
Then the generic fiber $Z_K$ contains precisely four $K$-rational points.
\item
Suppose that either $r=0,s\neq 0$ or that $\sqrt{r}\not\in K,s=0$
or that $r,s\neq 0$, $\frac{s}{r}\not\in\FF_4$, $\sqrt{r}\not\in K$.
Then $Z_K$ contains only one $K$-rational point.
\item
In all other cases, $Z_K$ contains exactly two $K$-rational points.
\end{enumerate}
\end{proposition}

\proof
Let $z\in Z_K$  be a $K$-rational point. Its preimage under the quotient map
$Y_L\ra Z_K$ is a closed subscheme of $K$-length two, which is invariant under
the diagonal $\alpha_2$-action. Hence the preimage must be an $\alpha_2$-invariant $L$-rational
point on $  Y_L$. The latter might be viewed as an equivariant section
for the structure map $Y_L\ra\Spec(L)$.
Since the  diagram 
$$
\begin{CD}
C_K @<<< Y_L\\
@VVV @VVV\\
\Spec(K) @<<< \Spec(L)
\end{CD}
$$
is cartesian, such sections are   equivariant $K$-morphisms
$h:\Spec(L)\ra C_K$. Note that the $\alpha_2$-actions on $C_K$ and $L$ are given by the
derivations $(u^{-2}+r)D_u$ and $(v^{-2}+s)D_v$, respectively.
Summing up, the $K$-rational points on $Z_K$ correspond to the equivariant
$K$-morphisms $h:\Spec(L)\ra C_K$.

Clearly, the image of $h$ must lie in the smooth locus of $C_K$.
So $h$ is given by a ring homomorphism
$$
K[u^{-1}]\lra L,\quad u^{-1}\longmapsto \alpha v^{-1}+\beta
$$
for some $\alpha,\beta\in K$.
Equivariance means that $(u^{-2}+r)D_u (u^{-1})=u^{-4}+ru^{-2}$
maps to $(v^{-2}+s)D_v(\alpha v^{-1}+\beta)=\alpha(v^{-4}+sv^{-2})$,
which gives the equation
$$
(\alpha v^{-1}+\beta)^4=\alpha(v^{-4}+sv^{-2})
$$
Comparing coefficients, we obtain three equations 
$$
\alpha^4=\alpha,\quad
r\alpha^2=\alpha s,\quad
\beta^2(\beta^2+r)=0.
$$
The first equation simply means $\alpha\in\FF_4$.
Unraveling the other two equations, we reach the assertion.
For example, consider the case that both $r,s\neq 0$.
The second equation then means $\alpha=0$ or $\alpha=\frac{s}{r}$.
This gives two or one possibilities for $\alpha$,
depending on whether $\frac{s}{r}$ is contained in $\FF_4$ or not.
The third equation means that $\beta=0$ or $\beta=\sqrt{r}$.
Again we have two or one possibility, depending on whether $r$ is
a square in $K$ or not.
Summing up, we have
either one, two, or four rational points in the generic fiber $Z_K$, depending on the values $\frac{s}{r}$ and $\sqrt{r}$ 
as claimed in the assertion.
The other cases are similar, and left to the reader.
\qed

\medskip
For example, in the case $r=s=1$, the homomorphism $K[u^{-1}]\ra L$ given
by $u^{-1}\mapsto v^{-1}$ corresponds to the rational point on $Z_K$
given by $a=b$, $c=0$, with maximal ideal   $\maxid=(a+b)$.

\begin{remark}
\mylabel{sections perfect}
For later use we record the outcome if  $\sqrt{r}\in K$, for example if $K$ is perfect.
Then the number of $K$-rational points in the generic fiber $Z_K$ depends only on
the ratio $(r:s)$ viewed as a point in the projective line $\PP^1$.
If $r=0$ we have precisely one rational point.
If $s=0$, or $r,s\neq 0$ and $(r:s)\not\in\PP^1(\FF_4)$ we have precisely two rational points.
And in case $r,s\neq 0$, $\frac{s}{r}\in\PP^1(\FF_4)$ we have four rational points.
\end{remark}

Next, we look at the closed fibers $Z_b=g^{-1}(b)$ of the quasielliptic fibration $g:Z\ra\PP^1$.
The situation is particularly simple for fibers over which the $\alpha_2$-action is free:

\begin{proposition}
\mylabel{multiplicity one}
Let $b\in\PP^1$ be a rational point so that $Y_b\subset Y$ contains no fixed points.
Then $Z_b\subset Z$ is isomorphic to the cuspidal curve of arithmetic genus one.
\end{proposition}

\proof
I make the computation for the special case $b=0$, the other cases
being similar.
Consider the commutative diagram (\ref{generic square}).
The fiber $h^{-1}(0)\subset C$ is the spectrum of the Artin ring
$A=k[v^2,v^3]/(v^2)$. In other words, we have $A=k[\epsilon]$,
where $\epsilon$ denotes the residue class of $v^3$. Our derivation $\delta=(u^{-2} +r)D_u + (v^{-2}+s) D_v$
acts as $\delta(\epsilon)=\epsilon$. The fiber $Y_A\subset Y$ is hence equivariantly isomorphic
to $C\times\alpha_2$ with diagonal $\alpha_2$-action.
We have  $Y_A/\alpha_2=Z_0$ because the group scheme action is free on $Y_A$, by
Proposition \ref{base free}.
Clearly, the quotient $Y_A/\alpha_2$ is the rational cuspidal curve of arithmetic genus one, hence
the assertion.
\qed

\medskip
Next, we come to multiple fibers:

\begin{proposition}
\mylabel{multiplicity two}
Let $b\in\PP^1$ be the rational point defined by the maximal ideal $(v^{-2})$.
Then the curve $Z_b$ is nonreduced, and we have an equality $Z_b=2(Z_b)_\red$ of Weil divisors. 
\end{proposition}

\proof
Let $h:C\ra\PP^1$ be the canonical morphism from Diagram (\ref{total square}).
The fiber $h^{-1}(b)\subset C$ is the spectrum of the Artin ring
$A=k[v^{-1}]/(v^{-2})$. In other words, we have $A=k[\epsilon]$,
where $\epsilon$ denotes the residue class of $v^{-1}$. Our derivation $\delta=(u^{-2} +r)D_u + (v^{-2}+s) D_v$
acts as $\delta(\epsilon)=0$. The fiber $Y_A\subset Y$ is therefore equivariantly isomorphic
to $C\otimes k[\epsilon]$, where $\alpha_2$ acts via the derivation $(u^{-2} +\lambda)D_u$ on the first factor, and
trivially on the second factor. Hence the quotient is isomorphic to $\PP^1\otimes k[\epsilon]$.
According to Proposition \ref{base free}, this quotient and the fiber  $Z_b$ are isomorphic on a dense open set,
and the assertion follows.
\qed

\medskip
In case that the parameter $s\in k$ is a square, the fiber $Z_b$ over the closed point
$b\in\PP^1$ defined by the maximal ideal $(v^{-2}-s)$ is  a double fiber as well.
This follows from Proposition \ref{automorphism group}.
The situation is different if $s\in k$ is not a square, which can happen only  over
nonperfect ground fields:

\begin{proposition}
Suppose that $s\in k$ is not a square.
Let $b\in \PP^1$ be the rational point defined by
the maximal ideal $(v^{-2}-s)$. Then the irreducible curve $Z_b$ is integral,
but not geometrically integral.
\end{proposition}

\proof
We just discussed that the fiber $Z_b$ has geometric multiplicity two.
The fiber $h^{-1}(b)\subset C$ is the spectrum of the Artin ring $k'=k[v^{-1}]/(v^{-2}-s)$, which is a field.
As in the preceding proof, we argue that the fiber $Z_b$ is birational to
the integral curve $\PP^1_{k'}$.
\qed

\medskip
We shall see in Section \ref{simultaneous} that the existence of such integral but geometrically multiple fibers
is responsible for the nonexistence of simultaneous resolutions.  The following result
also gives a hint that some purely inseparable base change is necessary:

\begin{proposition}
\mylabel{generic isomorphism}
Let $Z\ra\PP^1$ be the quasielliptic surface defined by parameters
$r,s\in k$, and $Z'\ra\PP^1$ be the quasielliptic surface defined
by $r'=t^2r, s'=t^2s$ for some nonzero $t\in k$.
Then there is a commutative diagram
$$
\begin{CD}
Z_\eta @>>> Z'_\eta\\
@VVV  @VVV\\
\Spec \k(\eta) @>>> \Spec \k(\eta),
\end{CD}
$$
such that the horizontal maps are isomorphisms, and $\eta\in\PP^1$ is
the generic point.
\end{proposition}

\proof
This is an exercise in Weierstrass equations.
The generic fibers $Z_\eta,Z'_\eta$ are isomorphic to 
regular cubics in $\PP^2_\eta$ containing a rational point.
According to Proposition \ref{quotient embedding}, the Weierstrass equation for $Z_\eta$
is
$$
y^2= (1+s^2b^2)x^3 + r^2b^3x^2 +b^3,
$$
where we set $y=c$ and $ x=a$, and $b\in k(\eta)$ is a transcendental generator.
Applying the  substitutions
$b\mapsto t^2b$  and $y\mapsto t^3y$ and $x\mapsto t^2x$, we obtain
the corresponding Weierstrass equation 
for $Z'_\eta$ up to a factor $t^6$.
\qed

\section{K3 surfaces and rational surfaces}
\mylabel{k3}

We keep the notation of the preceding two sections,
such that $Z$ is a proper normal surface, defined as the
quotient $Z=Y/\alpha_2$ and  depending on two parameters $r,s\in k$.
We also assume that our ground field $k$ is perfect.
Let $r:X\ra Z$ be the minimal resolution of singularities.
Note that we now have three surfaces $Y,Z,X$, and all three are
Cohen--Macaulay, hence each has a  dualizing sheaf. 

\begin{proposition}
\mylabel{gorenstein}
The proper normal surfaces $Y$ and $Z$ are Gorenstein, and their dualizing sheaves
are trivial as invertible sheaves.
If at least one parameter $r,s\in k$ is nonzero, the same holds for $X$.
\end{proposition}

\proof
The surface $Y=C\times C$ is Gorenstein and the invertible sheaf $\omega_Y$ trivial,
because the same holds for the curve $C$.
Let $U\subset Y$ be the complement of the fixed locus, and
$V\subset Z$ be the corresponding open subset.
The induced morphism $U\ra V$ is an
$\alpha_2$-torsor. By Proposition \ref{dualizing quotient} the quotient $V$ is Gorenstein and the
invertible sheaf $\omega_V$ is trivial.
Using that $Y$ is Cohen--Macaulay, we infer that the invertible sheaf $\omega_Y=\O_Y$ is 
trivial as well.

Suppose the  parameters $r,s\in k$ do not   vanish simultaneously. We saw
in Sections \ref{singularity} and \ref{singularities} that the singularities
on $Z$ are then rational double points, hence the relative dualizing sheaf $\omega_{X/Z}$
is trivial. We conclude that $\omega_X\simeq\O_X$.
\qed

\medskip
Using the Enriques classification of surfaces, it is now easy to determine 
the nature of the smooth proper surface $X$:

\begin{theorem}
\mylabel{k3 surface}
The smooth proper surface $X$ is a K3 surface if at least one parameter $r,s\in k$ is nonzero. 
Otherwise it is a geometrically rational surface.
\end{theorem}

\proof
We may assume that the ground field is algebraically closed.
Suppose $r,s\in k$ do not   vanish simultaneously. We saw in Proposition \ref{gorenstein}
that $\omega_X$ is trivial.
It follows that the regular surface $X$ contains no $(-1)$-curves, in other words,
$X$ is minimal.
By the Enriques  classification of surfaces, $X$ is either an abelian surface,
or an Enriques surface, or a K3 surface. Such surfaces have
second Betti number $b_2=6$, $b_2=10$, and $b_2=22$, respectively.
In Sections \ref{singularity} and \ref{singularities} we showed that the normal surface $Z$ contains
either five rational double points of type $D_4$, or one rational double
point of type $D_4$ and two of type $D_8$.
In any case, these singularities contribute twenty exceptional curves on $X$, and this implies
$b_2(X)\geq 20$. It follows that $X$ must be a K3 surface.

Now suppose $r=s=0$.
We shall apply the Castelnuovo Criterion for rationality.
According to Proposition \ref{elliptic singularity}, the normal surface
$Z$ contains an elliptic singularity.
Let $E\subset X$ be the corresponding exceptional divisor, and
$D\subset X$ be a curve supported by $E$. Using
$
K_{X/Z}\cdot D + D^2= K_X\cdot D + D^2 =\deg(\omega_D) = 2p_a(D)-2,
$
we infer $K_{X/Z}\cdot D\geq 0$. Since the intersection form 
on the divisors supported on $E$ is negative definite, the relative canonical
class $K_{X/Z}$, which is a divisor supported on $E$, has negative
coefficients. 
If we choose   $D$ with arithmetic genus $p_a(D)=1$, we have
$K_{X/Y}\cdot D = -D^2>0$. The upshot is that $K_X=K_{X/Z}$ is negative,
and all plurigenera $P_n(X)=h^0(X,\omega_X^{\otimes n})$, $n\geq 0$ vanish.
By Serre duality, $H^2(X,\O_X)$ vanishes, hence the Picard scheme $\Pic_{X/k}$ is reduced.
Since we have a quasielliptic fibration $X\ra \PP^1$, the Albanese map for $X$
is trivial, and we conclude that $H^1(X,\O_X)=0$. 
Now the Castelnuovo Criterion tells us that $X$ is rational.
\qed

\medskip
We are mainly interested in K3 surfaces. Let us therefore
assume that $r,s\in k$ do not   vanish simultaneously, such that $X$ is a K3 surface.
Let $f:X\ra\PP^1$ be the quasielliptic fibration induced by $g:Z\ra\PP^1$.
Suppose for the moment that the ground field $k$ is algebraically closed.
Obviously, quasielliptic  K3 surfaces are \emph{unirational}, that is, there is a surjective 
morphism from a rational surface onto $X$. 
In our situation, we may choose as rational surface the fiber product $X\times_Z \tilde{Y}$,
where $\tilde{Y}=\PP^1\times\PP^1$  is the normalization of $Y$.
Unirational K3 surfaces are \emph{supersingular} (in the sense of Shioda),
that is, the Picard number is $\rho(X)=22$.
There is a slight complication over nonclosed ground fields:

\begin{proposition}
\mylabel{picard number}
We have $\rho(X)=22$ if the ground field $k$ contains a third root
of unity, and $\rho(X)=21$ otherwise.
\end{proposition}

\proof
We    check this with the Tate--Shioda formula.
For each $b\in\PP^1$, let $\rho(X_b)$ be the number of reducible components
in the fiber $X_b\subset X$.
Clearly, $\sum(\rho(X_b) -1)$ is the number of irreducible components
in the exceptional divisor for the resolution of singularities $X\ra Z$.
It follows from the results in Sections \ref{singularity} and \ref{singularities}
that this sum 
equals 20 if the ground field $k$ contains a third root of unity, and 19 otherwise.
The group $\Pic^0(X_\eta)$ for the generic fiber $X_\eta$ is 2-torsion.
Hence the Tate--Shioda formula for the quasielliptic fibration $f:X\ra\PP^1$ takes the form
$$
\rho(X) = 2 + \sum(\rho(X_b) -1),
$$
and the result follows.
\qed

\medskip
Our next task is to understand the intersection form on the
N\'eron--Severi group $\NS(X)$, which for K3 surfaces is isomorphic to
the Picard group.
For this, we have to analyse the reducible fibers of the quasielliptic
fibration  $f:X\ra\PP^1$.

Recall that reducible fibers in genus-one fibrations with regular total
space $X$ correspond to root systems. The correspondence is as follows:
Decompose the fiber $X_b=E_0+E_1+\ldots +E_n$ into integral components,
and let $V$ be the real vector space generated by the $E_i$, endowed with
the positive semidefinite bilinear form $\Phi(E_i,E_j)=-(E_i\cdot E_j)$.
Suppose that $E_0$ is a component with multiplicity one.
Then the remaining $E_1,\ldots, E_n\in V$ form a root basis.
Note that $X_b-E_0$ is then the longest root, and its multiplicities
can be read off from the Bourbaki tables   \cite{Bourbaki LIE 4-6}.
One usually uses the
symbols for   positive semidefinite Coxeter matrix ($\tilde{A}_n$, $\tilde{B}_n, ...,\tilde{G}_2$)
to denote the fiber type.

In our situation it is very simple to determine  the fiber type of
$X_b$, $b\in\PP^1$.
 We only need to know the types of singularities $z\in Z$ lying in the fiber $Z_b$,
which we determined in Sections \ref{singularity} and \ref{singularities},
and the multiplicity of the closed fiber $Z_0\subset Z$, which we did
in Section \ref{quasielliptic}.
The proofs for the following two assertions are straightforward whence omitted.

\begin{proposition}
\mylabel{fiber d4}
Suppose $Z_b\subset Z$ is the fiber containing the rational double point $z\in Z$
corresponding to the quadruple point $y\in Y$ .
Then $X_b\subset X$ is a fiber of type $\tilde{D}_4$ if the ground field $k$ contains
a third root of unity, and of type $\tilde{B}_3$ otherwise.
The strict transform of the closed fiber $Z_b$ corresponds
to the white vertices in Figure 4.
\end{proposition}

\vspace{1em}
\centerline{\includegraphics{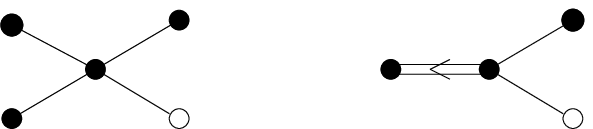}} 
\vspace{.8em}
\centerline{Figure \stepcounter{figure}\arabic{figure}: The extended Dynkin diagrams $\tilde{D}_4$ and $\tilde{B}_3$.}
\vspace{1em}

\medskip
We next examine fibers related to the fixed points of the $\alpha_2$-action on $X$.

\begin{proposition}
\mylabel{fiber d4d8}
Suppose $Z_b\subset Z$ is a fiber containing two rational double points $z \in Z$, say of 
type $D_4$ or $D_8$. 
Then $X_b\subset X$ is a fiber of type $\tilde{D}_8$ or $\tilde{D}_{16}$, respectively.
The strict transform of the closed fiber $Z_b$ corresponds
to the central white vertex, as in Figure 5.
\end{proposition}

\vspace{1em}
\centerline{\includegraphics{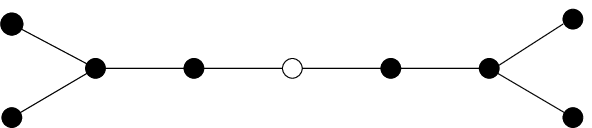}} 
\vspace{.8em}
\centerline{Figure \stepcounter{figure}\arabic{figure}: The extended Dynkin diagram $\tilde{D}_8$.}
\vspace{1em}

\begin{proposition}
\mylabel{fiber e8}
Suppose $Z_b\subset Z$ is a fiber containing precisely one  rational double point $z \in Z$  
corresponding to a fixed point  on  $Y$.
Then $X_b\subset X$ is a fiber of type $\tilde{E}_8$, and
the strict transform of the closed fiber $Z_b$ corresponds
to the white vertex in Figure 6.
\end{proposition}

\vspace{1em}
\centerline{\includegraphics{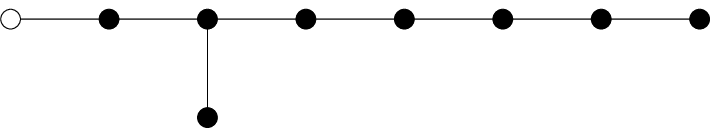}} 
\vspace{.8em}
\centerline{ Figure \stepcounter{figure}\arabic{figure}: The extended Dynkin diagram $\tilde{E}_8$.}
\vspace{1em}

\proof
Clearly, only the fiber types $\tilde{E}_8$ or $\tilde{D}_8$ are possible.
In the latter case, the strict transform of the closed fiber $Z_b$ must be
one of the four outer vertices. The outer vertices, however,
appear with multiplicity one in the fiber, contradicting that $Z_b$ has multiplicity two. 
\qed

\section{Discriminants and Artin invariants}
\mylabel{discriminants}

Let $X$ be a K3 surface over an algebraically closed ground field $k$ of characteristic $p>0$
with Picard number $\rho(X)=22$. In other words, $X$ is supersingular.
Artin \cite{Artin 1974a} introduced an integer invariant
for such surfaces called the \emph{Artin invariant}.
One way to define it is in terms of the intersection form on the
N\'eron--Severi group $\NS(X)$.
Artin showed that its discriminant is of the form
$\disc\NS(X)=-p^{2\sigma_0}$ for some integer  $1\leq\sigma_0\leq 10$,
which is the Artin invariant $\sigma_0(X)=\sigma_0$.
 
It is not difficult to compute the Artin invariant in the presence
of a quasielliptic fibration $f:X\ra\PP^1$ with a section $A\subset X$.
Let $L\subset\NS(X)$ be the subgroup generated by the section $A\subset X$,
together with all curves $C\subset X$ inside the closed fibers of 
the quasielliptic fibration $f:X\ra\PP^1$.
As Ito explained in \cite{Ito 1992}, Section 2, the quotient group
$\NS(X)/L$ is a finite group annihilated by $p$, say of order $p^n$.
This group acts freely and transitively on the set of rational points
in the generic fiber $X_\eta$. Hence the group order $p^n$ is
also the number of rational points on $X_\eta$.

For each point $b\in\PP^1$, let $C_b\subset X_b$ be the unique integral component
with $C_b\cdot A=1$. 
The remaining irreducible components $E_i\subset X_b-C_b$ form a root basis, whose type
corresponds to the fiber type of $X_b$.
Let $d_b=\det(E_i\cdot E_j)$ be the determinant of the corresponding intersection form.
 Bourbaki calls  this number  the \emph{connection index} (\cite{Bourbaki LIE 4-6}, chapter IV, \S 1.9).
Clearly, there are only finitely many points $b_1,\ldots, b_m\in\PP^1$ whose fibers
are reducible. To simplify notation, we set $d_i=d_{b_i}$.

\begin{lemma}
\mylabel{artin formula}
Under the preceding assumptions, the Artin invariant $\sigma_0 $ for the quasielliptic K3 surface $X$ is
given by the formula $p^{2\sigma_0+2n}=\pm d_1\ldots d_m.$
\end{lemma}

\proof
By definition,  the subgroup $L\subset\NS(X)$ has index $p^n$.
According to \cite{Serre 1979}, Chapter III, \S2, Proposition 5, we have   $\disc\NS(X) = \pm p^{-2n}\disc(L)$,
and it remains to compute the discriminant of $L$. 

Let $L'\subset L$ be the subgroup generated by the section $A\subset X$ and
a closed fiber $X_b\subset X$. Using $A^2=-2$ and $A\cdot X_b=1$ and $X_b\cdot X_b=0$,
we have $\disc(L')=-1$.
For each $1\leq i\leq m$ let $L_i\subset L$ be the subgroup generated by
the integral components of $X_{b_i}-C_{b_i}$, where $C_{b_i}\subset X_{b_i}$ is the
unique integral component with $C_{b_i}\cdot A=1$.
It is easy to see that we have an orthogonal decomposition
$L'\oplus L_1\oplus\ldots\oplus L_m$, and the result follows.
\qed

\medskip
We now return to the situation of the preceding section, such that $X$ is our
quasielliptic K3 surface in characteristic $p=2$ depending on two parameters $r,s\in k$,
which do not vanish simultaneously.
We then obtain a   point $(r:s)\in\PP^1(k)$.
The Artin invariant of $X$ now depends on whether or not this point
lies inside the subset $\PP^1(\FF_4)\subset\PP^1(k)$.

\begin{proposition}
\mylabel{}
Our supersingular K3 surface $X$ has Artin invariant
$\sigma_0 =1$ if $(r:s)\in\PP^1(\FF_4)$.
Otherwise, it has Artin invariant $\sigma_0=2$.
\end{proposition}

\proof
Suppose first that $r=0$.
According to Proposition \ref{sections} and Remark \ref{sections perfect}, the generic fiber $X_\eta$
contains only $1=2^0$ rational point.
By Propositions \ref{fiber d4} and \ref{fiber e8}, there are
three reducible fibers, one of type $\tilde{D}_4$, the others of type $\tilde{E}_8$.
The connection indices are $d_1=-4$ and $d_2=d_3=1$.
By the formula in Lemma \ref{artin formula}, we have $\sigma_0=1$.

Now suppose that $r\neq 0$, $s=0$. Then the the generic fiber
contains precisely $2=2^1$ rational points. By Proposition \ref{},
we have one fiber of type $\tilde{D}_4$ and one fiber of type $\tilde{D}_8$, which both have connection index
$d_1=d_2=-4$. The formula for the Artin invariant again gives
$\sigma_0=1$. 

In case that $r,s\neq 0$ but $\frac{s}{r}\not\in\FF_4$ we also have only two rational points
in the generic fiber, but three fibers of type $\tilde{D}_4$.
This implies $\sigma_0=2$.

Finally, suppose we have $r,s\neq 0$ and $\frac{s}{r}\in\FF_4$.
Then we have $4=2^2$ rational points in the generic fiber and again three fibers
of type $\tilde{D}_4$, and this leads to  $\sigma_0=1$.
\qed

\begin{remark}
There seems to be a close relation between our family of K3 surfaces and the
family studied in \cite{Schroeer 2004}.
\end{remark}

\section{Blowing up curves on rational singularities}
\mylabel{blowing}

The next task is to construct   simultaneous resolutions of singularities
for our flat family of normal K3 surfaces.
By the work of Brieskorn \cite{Brieskorn 1968} and Artin \cite{Artin 1974b},
simultaneous resolutions in flat families rarely exist without base change.
In our case, it turns out that a purely inseparable base change is necessary.
After that, simultaneous resolution is achieved by blowing up Weil divisors inside the quasielliptic fibration,
which easily extends to the   family. The goal of this section
is to collect some useful facts on blowing up curves on rational surface singularities.
Throughout, $(\O_z,\maxid,k)$ is a 2-dimensional normal local ring, with 
resolution of singularities $r:X\ra\Spec(\O_z)$. 
We  assume that $\O_z$ is a \emph{rational singularity}, that is, $H^1(X,\O_X)=0$.

Let $I\subset\O_z$ be a reflexive ideal, which defines a curve $C\subset\Spec(\O_z)$
without embedded components. How to compute the schematic fiber
$r^{-1}(C)\subset S$ on the resolution of singularities?
The following arguments are adapted from Artin's paper \cite{Artin 1966}, where
he considered maximal ideals instead of reflexive ideals.

Let $E\subset X$ be the reduced exceptional divisor, and $E=E_1+\ldots+E_n$
be its decomposition into integral components, and $C'\subset X$ be the
strict transform of $C\subset\Spec(\O_z)$. Consider nonzero divisors $Z=\sum r_i E_i$
with coefficients $r_i\geq 0$ satisfying $(Z+C')\cdot E_i\leq 0$ for all $1\leq i\leq n$.
As in \cite{Artin 1966}, page 131 there is a unique minimal cycle $Z$ with these properties.
(In Artin's situation, this cycle is called the \emph{fundamental cycle}.)

One may determine this cycle by computing  a sequence of cycles $Z_0,Z_1,\ldots$ inductively
 as follows: Start with $Z_0=E_i$, where $E_i$ is any component with $E_i\cdot C'>0$. Suppose we already defined $Z_m$ for some $m\geq 0$.
If $(Z_m+C')\cdot E_i\leq 0$ for all $1\leq i\leq n$, we set $Z=Z_m$ and are done.
Otherwise, we have $(Z_m+C')\cdot E_i>0$ for some integral component $E_i$.
We then define $Z_{m+1}=Z_m+E_i$ and proceed by induction. This algorithm stops after
finitely many steps and yields the desired cycle $Z$.

\begin{lemma}
\mylabel{fundamental cycle}
We have $r^{-1}(C)=Z\cup C'$ as subschemes of the resolution of singularities $X$.
\end{lemma}

\proof
The arguments are as in the proof for \cite{Artin 1966}, Theorem 4.
\qed

\begin{remark}
This result shows that in the algorithm to compute $Z$, we may start by letting $Z_0$
be the fundamental cycle of the singularity $\O_z$. Note that the fundamental cycle
corresponds to the longest roots, which can be read off from the Bourbaki tables \cite{Bourbaki LIE 4-6}.
\end{remark}

\begin{remark}
Mumford \cite{Mumford 1961} defined a linear preimage from the group of
1-cycles on $\Spec(\O_z)$ to the group of 1-cycles with rational coefficients on $S$.
Note that the schematic fiber $r^{-1}(C)\subset X$ usually differs from the
linear preimage $r^*(C)\in\Div(X)\otimes\QQ$.
\end{remark}

\begin{lemma}
\mylabel{two generators}
Let $I$ be a reflexive fractional $\O_z$-ideal. Then the $k$-vector space $I\otimes_{\O_z} k$
is at most 2-dimensional.
\end{lemma}

\proof
We may assume that the local ring $\O_z$ is henselian, and that the fractional
ideal $I$ is an ideal in $\O_z$,  which defines a local curve $C\subset\Spec(\O_z)$.
It induces an ideal sheaf $\shI=I\O_S$ on the resolution of singularities $S$.
I claim that the canonical map $I\ra H^0(S,\shI)$ is bijective. To see this,
consider the commutative diagram
$$
\begin{CD}
0 @>>> I @>>> \O_z @>>> \O_z/I @>>> 0\\
@. @VVV @VVV @VVV\\
0 @>>> H^0(X,\shI) @>>> H^0(X,\O_X) @>>> H^0(X,\O_X/\shI).
\end{CD}
$$
The vertical map in the middle is surjective, because $\O_z$ is normal and $X\ra\Spec(\O_z)$
is proper and birational. The vertical map on the right is injective, because $\O_z/I$
has no embedded components. Now the Snake Lemma tells us that $I\ra H^0(X,\shI)$
is bijective.

I claim that there are two global sections $f,f'\in H^0(X,\shI)$ generating the stalk $\shI_x$
for all points $x\in X$.  To construct $f$, choose a divisor $D_0\subset E$ whose support
is disjoint from $x$ and $\Sing(E)$, and whose degree on $E_i$ equals the intersection number
$-r^{-1}(C)\cdot E_i$ for all $1\leq i\leq n$. Note that the closed subscheme $r^{-1}(C)\subset X$ indeed has
no embedded points, hence is  a Cartier divisor, according to Lemma \ref{fundamental cycle}. Since $\O_z$ is henselian, we may extend
$D_0\subset E$ to a Cartier divisor $D\subset S$ by
\cite{EGA IVd}, Proposition 21.9.11. By construction
$$
-D\cdot E_i = r^{-1}(C)\cdot E_i.
$$
Using that $\O_z$ is rational, we conclude that the invertible sheaf $\O_X(D + r^{-1}(C)) $ is trivial,
whence there is an isomorphism $f:\O_X(D)\ra\shI$.
Repeating the same construction, we find another divisor $D'\subset X$ as above, and whose
support is disjoint from $D$. This yields the desired surjection $f+f':\O_X^{\oplus 2}\ra\shI$.

We obtain an exact sequence $0\ra \shL\ra\O_X^{\oplus 2}\ra\shI\ra 0$,
which induces a long exact sequence
$$
H^0(X,\O_X^{\oplus 2})\lra H^0(X,\shI)\lra H^1(X,\shL).
$$
The sheaf $\shL$ is invertible. Taking determinants, one sees that
$\shL$ is isomorphic to the dual $\shI^\vee=\O_X(r^{-1}(C))$.
Consequently, $\shL\cdot E_i=r^{-1}(C)\cdot E_i\leq 0$.
Using \cite{Giraud 1982}, Proposition 1.9, we infer $H^1(X,\shL)=0$.
Summing up, $H^0(X,\shI)$ and hence $I$ are generated by two elements.
\qed

\medskip
Next, consider  the class group  $\Cl(\O_z)$ of reflexive fractional ideals $I$.
Given a class $[I]\in\Cl(\O_z)$, we form the blowing up 
$$
Z'=\Proj(\bigoplus_{n\geq 0} I^n).
$$
It depends, up to isomorphism, only on the ideal class and not the ideal itself.
Note that the induced map $h:Z'\ra\Spec(\O_z)$ is projective and birational.

\begin{proposition}
\mylabel{factorization blowing}
The resolution of singularities $r:X\ra\Spec(\O_z)$ factors
over our  blowing up $h:Z'\ra \Spec(\O_z)$.
\end{proposition}

\proof
We may assume that the fractional ideal $I$ is an ideal in $\O_z$,
defining a curve $C\subset\Spec(\O_z)$. According to Lemma \ref{fundamental cycle},
the schematic preimage $r^{-1}(C)\subset X$ is a Cartier divisor.
By the universal property of blowing ups, this means that $r$ factors over $h$.
\qed

\medskip
The next result tells us that the induced morphism $X\ra Z'$
coincides with its Stein factorization, whence is uniquely determined
by its exceptional curves.

\begin{proposition}
\mylabel{normal blowing}
The scheme $Z'$ is normal.
\end{proposition}

\proof
We may assume that $I$ is an ideal in $  \O_z$. The quotient $\O_z/I$ has no
embedded primes. 
Let $Z''\ra\Spec(\O_z)$ be any proper birational morphism, with $Z''$ reduced.
In light of Lipman's results \cite{Lipman 1969},
it suffices to check that the canonical map $I\ra\Gamma(Z'',I\O_{Z''})$
is bijective. For this, we argue as in Lemma \ref{two generators}.
\qed

\medskip 
We infer that the fibers of the blowing up are as small as possible:

\begin{proposition}
\mylabel{blowing fiber}
Suppose that $[I]\neq 0$. Then the closed
fiber $h^{-1}(z)\subset Z'$ is isomorphic to the projective line $\PP^1_k$.
\end{proposition}

\proof
Write $I=(f,g)$. This defines a closed embedding $Z'\subset\PP^1\times\Spec(\O_z)$,
hence the fiber $h^{-1}(z)\subset Z'$ is a closed subset of $\PP^1$.
Suppose we have $h^{-1}(z)\neq \PP^1$. Then the closed fiber is finite,
hence $h:Z'\ra\Spec(\O_z)$ is finite and birational.
Since $\O_z$ is normal, we conclude that $Z'=\Spec(\O_z)$.
In turn, $I$ must be invertible, contradiction.
\qed

\medskip
In light of the preceding results, it is possible to compute
the blowing up $Z'\ra \Spec(\O_z)$ of some reflexive fractional ideal $I$
via a resolution of singularities $r:X\ra \Spec(\O_z)$ as follows:
Let $E\subset X$ be the reduced exceptional divisor, and assume that $I$ is an ideal
defining a curve $C\subset\Spec(\O_z)$. The scheme $Z'=\Proj(\bigoplus_{n\geq 0} I^n)$
is obtained from $X$ by contracting all integral components $E_i\subset E$ but one.
The next result tells us which ones:

\begin{proposition}
\mylabel{blown down}
The integral components $E_i\subset E$ that are contracted by
$X\ra Z'$  are precisely those with $r^{-1}(C)\cdot E_i=0$.
\end{proposition}

\proof
The Cartier divisor $h^{-1}(C)\subset Z'$ is $h$-antiample.
The projection formula implies that $r^{-1}(C)\cdot E_i<0$ if $E_i$
is not contracted, and $r^{-1}(C)\cdot E_i=0$ if $E_i$ is contracted. 
\qed

\section{Blowing ups in genus-one fibrations}
\mylabel{genus}

We now apply the results from the previous section to the following situation.
Suppose $X$ is a smooth surface endowed with a genus-one fibration
$f:X\ra\PP^1$. We assume that the fibration is relatively minimal,
admits a section, and that the generic fiber $X_\eta$ is  a regular
curve of arithmetic genus one. Fix a rational point $b\in\PP^1$, and
decompose the closed fiber $(X_b)_\red=E_0+\ldots+ E_n$ into integral
components. Let $r:X\ra Z$ be the contraction of some integral components in 
$X_b\subset X$. We seek to recover the minimal
resolution of singularities with a sequence 
$$
X=Z^{(n)}\lra \ldots \lra Z^{(1)}\lra Z^{(0)}=Z
$$
of blowing ups so that each step  $Z^{(i+1)}\ra Z^{(i)}$
has as center a Weil divisor $C^{(i)}\subset Z^{(i)}_b$ inside the fiber.
This approach has its merits when it comes to deformations: Such 1-dimensional
centers behave much better in families than 0-dimensional centers.
A convenient way to describe the centers $C^{(i)}$ is via its strict transforms on $X$.

Suppose now that the fiber $X_b$ has fiber type $\tilde{D}_8$,
and let $r:X\ra Z$ be the contraction of all integral components but $E_4\subset X$,
which corresponds to the white vertex as depicted in Figure 7.
As usual, the choice of indices is taken from the Bourbaki tables   \cite{Bourbaki LIE 4-6}.
Throughout, the integral component $E_i\subset X_b$ shall correspond to the vertex with number $i$.
Note that the situation is as in Proposition \ref{fiber d4d8}.

\vspace{1em}
\centerline{\includegraphics{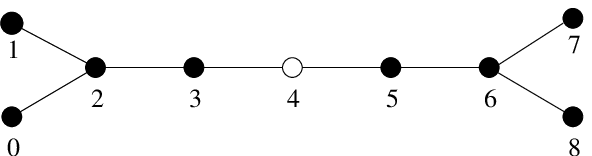}} 
\vspace{.8em}
\centerline{ Figure \stepcounter{figure}\arabic{figure}: The fiber type $\tilde{D}_8$.}
\vspace{1em}

\begin{proposition}
\mylabel{reach d8}
Under the preceding assumptions, we reach the   resolution of singularities
with the sequence of blowing ups $X=Z^{(6)}\ra \ldots \ra Z^{(0)}=Z$, in which
the centers $C^{(0)},\ldots, C^{(5)}$ have strict transforms
$E_4, 2E_5, E_3, E_2, E_5, E_6$, respectively.
All centers have arithmetic genus $p_a=0$.
\end{proposition}

\proof
By assumption, the normal surface $Z$ contains precisely two singularities, which are
rational double points of type $D_4$. Let $Z^{(1)}\ra Z$ be the blowing up whose center
$C^{(0)}$ has strict transform $E_4$. Note that this center is nothing but the half fiber.
Using the algorithm in Section \ref{blowing}, we compute its schematic
preimage $F^{(0)}\subset X$ on the minimal resolution of singularities.
It turns out that
$$
F^{(0)}=E_0+E_1+2E_2+2E_3 + E_4 + 2E_5+2E_6+E_7+E_8,
$$
which has $F^{(0)}\cdot E_3=F^{(0)}\cdot E_5=-1$.
Therefore, the exceptional curve for  $Z^{(1)}\ra Z$ corresponds to  $E_3+E_5$.
We now have to repeat this, quite mechanically. My findings are summarized in the following table:
$$
\begin{array}[t]{l |l| l |l}
        &\text{RDP} & \text{center}& \text{schematic preimage of center}\\ \hline&&&\\[-2ex]
Z^{(0)} & 2D_4       & E_4          & E_0+E_1+2E_2+2\underline{E_3} + E_4 + 2\underline{E_5}+2E_6+E_7+E_8  \\
Z^{(1)} & 2A_3       & 2E_5         & 2E_5+2E_6+E_7+E_8                                                         \\
Z^{(2)} & 2A_3       & E_3          & E_0+E_1+2\underline{E_2}+ E_3                   \\
Z^{(3)} & 2A_1+A_3   & E_2          & \underline{E_0}+\underline{E_1}+E_2\\
Z^{(4)} & A_3        & E_5          & E_5+2\underline{E_6}+ E_7+E_8\\
Z^{(5)} & 2A_1       & E_6          & E_6+\underline{E_7}+\underline{E_8}
\end{array}  
$$
\vspace{3em}
\centerline{Table \stepcounter{table}\arabic{table}: \text{Blowing ups in the $\tilde{D}_8$-fiber.}}
\vspace{-3.5em}

The second column displays the types of rational double points on $Z^{(i)}$. 
The third column gives the strict transforms  of the centers $C^{(i)}\subset Z^{(i)}$ on $X$.
The last column contains the schematic preimage  of the center $C^{(i)}$ on $X$.
The underlined irreducible components  are those with nonzero intersection number with the preimage, hence
give the exceptional curves for $Z^{(i+1)}\ra Z^{(i)}$. 

Note that the Weil divisor $C^{(1)}$ is already   Cartier, such that
$Z^{(2)}\ra Z^{(1)}$ is the identity. I have included this seemingly superfluous step here
because we need it later when it comes to simultaneous resolutions in families.

It remains to check $p_a(C^{(i)})=0$. I do this for $i=0$ and $i=1$, the other cases being similar.
The strict transform $E_4\subset X$ for $C^{(0)}\subset Z$ is clearly isomorphic to $\PP^1$, so the possibly nonnormal
points on $C^{(0)}$ must appear at the singularities $z\in Z$. To see that
$C^{(0)}$ is normal, it suffices to check that $E_4\cap r^{-1}(z)=\Spec(k)$.
According to \cite{Artin 1966}, Theorem 4, the fiber $r^{-1}(z)$ is the fundamental cycle
for the singularity,  which in our case equals $E_5+2E_6+E_7+E_8$, and
we see $E_4\cdot(E_5+2E_6+E_7+E_8)=1$.
It follows   $p_a(C^{(0)})=0$.

As above, one checks that  $C^{(1)}_\red\simeq\PP^1$.
The center $C^{(1)}$ and its strict transform $2E_5$ are nonreduced.
They are \emph{ribbons} in the terminology of Bayer and Eisenbud, that
is, infinitesimal extension of a reduced scheme by an invertible sheaf.
The strict transform $2E_5$ is the infinitesimal extension of $\PP^1$ by $\O_{\PP^1}(2)$.
Let $z\in Z^{(1)}$ be the singularity lying on $C^{(1)}$, which is a rational double point
of type $A_3$. As above, one shows that its preimage on $E_5$ is the intersection
$E_5\cap (2E_6+E_7+E_8)=\Spec(k[\epsilon])$.
This implies that $C^{(1)}$ is an extension of $\PP^1$ by $\O_{\PP^1}=\O_{\PP^1}(2-2)$,
as explained in \cite{Bayer; Eisenbud 1995}, Corollary 1.10.
The exact sequence
$
H^1(\PP^1,\O_{\PP^1})\ra H^1(C^{(1)},\O_{C^{(1)}})\ra H^1(\PP^1,\O_{\PP^1})
$
gives $p_a(C^{(1)})=0$.
\qed

\medskip
We now turn to the following case: Suppose that $X_b$ has fiber type $\tilde{E}_8$,
and that $X\ra Z$ is the contraction of all irreducible components $E_i\subset X_b$ except
$E_1$, which corresponds to the white vertex in Figure 8.
Note that the situation is precisely as in Proposition \ref{fiber e8}.

\vspace{1em}
\centerline{\includegraphics{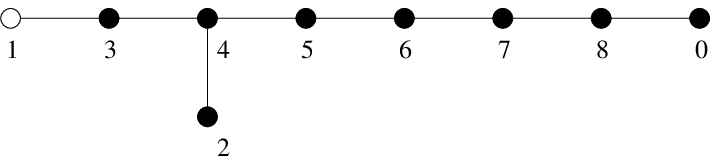}} 
\nopagebreak\vspace{.8em}\nopagebreak
\centerline{Figure \stepcounter{figure}\arabic{figure}: The fiber type $\tilde{E}_8$.}
\vspace{1em}

\begin{proposition}
\mylabel{reach e8}
Under the preceding assumptions, we reach the   resolution of singularities
with the sequence of blowing ups $X=Z^{(6)}\ra \ldots \ra Z^{(0)}=Z$, in which
the centers $C^{(0)},\ldots, C^{(5)}$ have strict transforms
$E_1, 2E_3, E_3, E_4, 2E_2+E_3+2E_4+2E_6,E_8$, respectively. All centers have arithmetic genus $p_a=0$.
\end{proposition}

\proof
By assumption, the normal surface $Z$ contains precisely one singularity, which is
a rational double point of type $D_8$. Let $Z^{(1)}\ra Z$ be the blowing up whose center
$C^{(0)}$ has strict transform $E_1$. Note that this is nothing but the half fiber.
Using the algorithm in Section \ref{blowing}, we compute its schematic
preimage $F^{(0)}\subset X$ on the minimal resolution of singularities.
It turns out that
$$
F^{(0)}= E_1+3E_2+4E_3+6E_4+5E_5+4E_6+3E_7+2E_8+E_0,
$$
which has $F^{(0)}\cdot E_3= -1$.
Therefore, the exceptional curve for  $Z^{(1)}\ra Z$ corresponds to  $E_3 $.
Further calculations are summarized as above in the following table:
$$
\hspace{-5em}
\begin{array}[t]{l |l| l |l}
        & \text{RDP} & \text{center}& \text{schematic preimage of center}\\ \hline&&&\\[-2ex]
Z^{(0)} & D_8         & E_1          & E_1+3E_2+4\underline{E_3}+6E_4+5E_5+4E_6+3E_7+2E_8+E_0 \\
Z^{(1)} & A_7         & 2E_3         & 2E_2+2E_3+4E_4+4E_5+4\underline{E_6} +3E_7+2E_8+E_0   \\
Z^{(2)} & 2A_3        & E_3          & E_2+E_3+2\underline{E_4}+E_5\\               
Z^{(3)} & 2A_1+A_3    & E_4          & \underline{E_2}+E_4+\underline{E_5}    \\
Z^{(4)} & A_3         & E_2+E_3+2E_4+2E_5+2E_6     &  E_2+E_3+2E_4+2E_5+2E_6+2E_7+2\underline{E_8}+E_0\\
Z^{(5)} & 2A_1        & E_8          & \underline{E_7}+E_8+\underline{E_0}
\end{array}
$$
\vspace{2em}
\centerline{Table \stepcounter{table}\arabic{table}:  Blowing ups in the $\tilde{E}_8$-fiber.}
\vspace{-2.5em}

The first assertion follows. Checking that $p_a(C^{(i)})=0 $  
is as in the proof for Proposition \ref{reach d8}.
\qed

\section{Simultaneous resolutions and nonseparatedness}
\mylabel{simultaneous}

Let $S\subset\Spec k[r,s]$ be the complement of the origin,
and $\foZ\ra S$ be our projective flat family of normal K3 surfaces with rational double points,
as constructed in Section \ref{selfproduct}.
In this section we seek a simultaneous resolution of singularities.
By Brieskorn's work, such simultaneous resolutions seldomly exist
without base change \cite{Brieskorn 1968}. Artin introduced the
\emph{resolution functor} $\Res_{\foZ/S}$ to remedy
the situation \cite{Artin 1974b}. Given an $S$-scheme $S'$, the resolution functor takes as values $\Res_{\foZ/S}(S')$
the set of isomorphism classes of simultaneous resolutions over $S'$, that is, commutative diagrams
$$
\begin{CD}
\foX' @>>> \foZ\\
@VVV @VVV\\
S' @>>> S,
\end{CD}
$$
where $\foX'\ra S'$ is a proper smooth family of K3 surfaces, and for each point $\sigma\in S'$,
the canonical morphism $\foX'_\sigma\ra\foZ_\sigma$ is the minimal resolution of singularities.
Artin's   insight was that the resolution functor is representable by an  algebraic space
that, however, is   only locally quasiseparated. The goal of this section is to show that
simultaneous resolutions exist after purely inseparable base change:

\begin{theorem}
\mylabel{simultaneous resolution}
There exist a simultaneous resolution $\foX'\ra S'$  for the family $\foX\ra S$ over
the purely inseparable flat base change
$S'=S\otimes_{k[r,s]} k[\sqrt{r},\sqrt{s}]$.
\end{theorem}

\proof
Recall that the quadruple point $y\in Y=C\times C$ induces a family of
$D_4$-singularities in $\foZ\ra S$. According to Proposition \ref{trivial deformation},
this family of singularities is formally trivial. It is then easy to see
that a simultaneous resolution already exists over $S$ without base change.

The trouble comes from the other singularities, which indeed necessitate a base change.
To proceed we cover the quasiaffine scheme $S'\subset\Spec k[\sqrt{r},\sqrt{s}]$
by the two affine open subsets
$U=D(\sqrt{s})$ and $V=D(\sqrt{r})$. Let us first concentrate 
on the restriction $\foZ_U\ra U$ to the first affine open subset.
As discussed in Section \ref{quasielliptic}, the second projection
$\pr_2:Y=C\times C\ra C$ induces a family of quasielliptic structures
$f:\foZ_U\ra \PP^1_U$. In accordance with our notation, we write the base
as $\PP^1_U=\Proj\O_U[v^2]$.
For each $\sigma\in U$, the induced fibration $f:\foZ_\sigma\ra\PP^1_\sigma$ has
over $v^{-2}=0$ and $v^{-2}=s$ degenerate fibers, of type $\tilde{D}_8$ over the open
subset $U\cap V$, and of type $\tilde{E}_8$ over the closed subset $U-V$.
It is precisely this point where   we use  base change, in order to apply Proposition \ref{automorphism group}.

Now let $\eta\in U$ be the generic point. We discussed in the preceding section
how to obtain the minimal resolution of $\foZ_\eta$ via a sequence of blowing ups
$$
\foX_\eta =\foZ_\eta^{(6)}\lra \ldots \lra\foZ_\eta^{(0)} = \foZ_\eta
$$
whose centers $\foC_\eta^{(i)}\subset\foZ_\eta^{(i)}$ are Weil divisors in the reducible fibers.
I claim that this procedure extends to the family $\foZ\ra U$.
Indeed: The first center $\foC_\eta=\foC_\eta^{(0)}$ is nothing but the half fiber
inside the degenerate fibers. Let 
$\foC_U\subset\foZ_U$   be the schematic closure of $\foC_\eta$.
On the smooth part of $\foZ_U\ra U$ this is a relative Cartier divisor,
and its restrictions $\foC_\sigma$, $\sigma\in U$ gives the half fibers.
At the   singularities, the schematic fibers $\foC_\sigma$ could pick up   embedded components.
This, however, is impossible, because $h^1(\O_{\foC_\sigma})=0$ by Propositions
\ref{reach d8} and \ref{reach e8}, and the Euler characteristic $\chi(\O_{\foC_\sigma})$ is constant by flatness.
We deduce that the blowing up of $\foC_U\subset\foZ_U$ yields fiberwise the first
step in the resolution of singularities. 

Repeating the preceding argument  inductively for the other centers,
we see that after six steps we reach a simultaneous resolution $\foX_U\ra\foZ_U$,
at least over the open subset $U\cap V$ where the fiber type is constantly $\tilde{D}_8$.
Using the tables in the preceding section  we can check that centers in the $\tilde{D}_8$-fibers specialize
to the centers in the $\tilde{E}_8$-fibers, hence we obtain the desired resolution of
singularities over $U$. Below is a table describing 
how the integral components of the $\tilde{D}_8$-fibers specialize to curves on 
the $\tilde{E}_8$-fiber. Note that this specialization  respects intersection numbers.
The existence of such a specialization seems to be interesting in its own right.
$$
\hspace{-1em}
\begin{array}[t]{l|l|l|l|l|l|c|l|l|l}
\text{$\tilde{D}_8$-fiber} & E_0 & E_1 & E_2 & E_3 & E_4  & E_5                       & E_6 & E_7 & E_8\\ 
\hline&&&\\[-2ex]
\text{$\tilde{E}_8$-fiber} & E_2 & E_5 & E_4 & E_3 & E_1 & E_2+E_3+2E_4+2E_5+2E_6+E_7 & E_8 & E_7 & E_0
\end{array}
$$
\vspace{2em}
\centerline{Table \stepcounter{table}\arabic{table}:  Specialization from $\tilde{D}_8$-fiber to $\tilde{E}_8$-fiber.}
\vspace{-2.5em}

The reason for this specialization behavior is as follows:
On $\foZ^{(0)}$, the $E_4$-curve in the $\tilde{D}_8$-fibers specialize to
the $E_1$-curve in the $\tilde{E}_8$-fibers.
We then do the first blowing up $\foZ^{(1)}\ra\foZ^{(0)}$. 
According to Tables 1 and Table 2 in Section \ref{genus}, its center is the
family of curves comprising the  $E_4$-curve on the $\tilde{D}_8$-fibers and the  $E_1$-curve on the $\tilde{E}_8$-fibers.
The preimage on $\foZ^{(1)}$ of the center 
is the family of curves  comprising  the $2E_3+E_4+2E_5$-curve on the $\tilde{D}_8$-fibers and
$E_1+4E_3$-curve on the $\tilde{E}_8$-fibers.
Consequently, both the $E_3$-component and the $E_5$-component specialize to the $E_3$-component.
The $E_4$-components specialize to the $E_1$-component.
Note that the preimage of the center is the Cartier divisor corresponding
to the canonical section into $\O_{\foZ^{(1)}}(-1)$.

Next, we do the second blowing up $\foZ^{(2)}\ra\foZ^{(1)}$. 
Now the center is the family comprising the $2E_5$-curve on the $\tilde{D}_8$-fibers  and
the $2E_3$-curve on the $\tilde{E}_8$-fibers.
The preimage on $\foZ^{(2)}$ of the center  is the family consisting of  the $2E_5$-curve on the $\tilde{D}_8$-fibers and
$2E_3+4E_6$-curve on the $\tilde{E}_8$-fibers.
Hence the $E_5$-component specializes to the $E_3+2E_6$-curve.
The $E_3$-components specializes to the $E_3$-component, and the $E_4$-component to the $E_1$-component.

Now we come to the third blowing up $\foZ^{(3)}\ra\foZ^{(2)}$. 
Here the center is the family comprising the $E_3$-curve on both the $\tilde{D}_8$-fibers  and
the  $\tilde{E}_8$-fibers. The preimage of the center is the family comprising
the $E_3+2E_2$-curves on the $\tilde{D}_8$-fibers and
$E_3+2E_4$-curve on the $\tilde{E}_8$-fibers. So the $E_3$-component specializes to $E_3$, and the  $E_2$-component specialize to
$E_4$. The $E_4$-component keeps on specializing to the $E_1$-component, because this family is
disjoint from the exceptional locus of the blowing-up.
The specialization of the $E_5$-component on the $\tilde{D}_8$-fibers is more interesting:
It is necessarily of the form
$E_3+nE_4+2E_6$ for some integer $n\geq 0$.
Using the fact that the whole fibers specializes to a whole fiber, we infer $n=2$.
The arguments for the remaining two steps in the blowing up process are similar, and left to the reader.

It remains to extend this simultaneous resolution  to $S'$.
Using the first projection $\pr_1:Y=C\times C\ra C$ rather than the second projection,
we obtain a simultaneous resolution $\foX_V\ra V$ for the family  $\foZ_V\ra V$.
The two resolutions $\foX_U\ra U$ and $\foX_V\ra V$ coincide over the overlap
$U\cap V$. Indeed, the degenerate fibers in $\foZ\ra\PP^1_S$ have constant
fiber type $\tilde{D}_8$ over $U\cap V$, and the fiberwise integral components of the  exceptional divisors
in the simultaneous resolution $\foX_U$ constitute relative Cartier divisors over $U\cap V$.
Such simultaneous resolutions are necessarily unique.
\qed

\medskip
The   result may be rephrased by saying that there is a morphism $S'\ra\Res_{\foZ/S}$. 
This does not seem to be an isomorphism, because our simultaneous resolution $\foX'\ra S'$ does not appear to
be unique. This nonuniqueness leads to nonseparability phenomena, as discussed
by Burns and Rapoport in \cite{Burns; Rapoport 1975}, Section 7.

\section{Isomorphic fibers in the family}
\mylabel{isomorphic}

In the preceding section we constructed a smooth family of supersingular K3-surfaces $\foX'\ra S'$
with Artin invariants $\sigma_0\leq 2$, which is defined
over the complement $S'\subset\AA^2$ of the origin, where $\AA^2=\Spec k[\sqrt{r},\sqrt{s}]$.
According to Rudakov and Shafarevich \cite{Rudakov; Safarevic 1978}, such families
should depend only on one effective parameter.
This is indeed the case, and can be made explicit as follows. Write $\PP^1=\Proj k[\sqrt{r},\sqrt{s}]$,
and consider the canonical projection $S'\ra\PP^1$.

\begin{theorem}
\mylabel{induced family}
Our smooth family $\foX'\ra S'$  is isomorphic to the pullback
of a smooth family $\foX\ra\PP^1$ of supersingular K3 surfaces with
Artin invariants $\sigma_0\leq 2$.
\end{theorem}

\proof
The assertion would be obvious if the  moduli space of polarized K3 surfaces
would be fine. This, however, does not seem to be the case.
Instead we shall work with moduli space of marked K3 surfaces.
For this we have to check that our family admits a marking.

Let $K=k(\sqrt{r},\sqrt{s})$ be the function field of the pointed
affine plane $S'$, and choose a separable closure $K\subset K^\sep$.
The Galois group $G=\Gal(K^\sep/K)$ acts on the the module $N=\Pic(X_{K^\sep})$,
where $X_K=\foX'_K$ denotes the generic fiber.
This action must be trivial: According to the explicit description
of   fibers and sections in $X_K$ in Sections \ref{singularity} and \ref{singularities}
and Proposition \ref{sections}, the pullback map
$\Pic( X_K)\ra\Pic( X_{K^\sep})$ is bijective.
We conclude that the sheaf $\Pic_{\foX'/S'}$ on $S'$ in the \'etale topology
is constant on some open affine subset $U\subset S'$.
Whence there is bijection $N_U\ra(\Pic_{\foX'/S'})_U$.
This bijection extends to a homomorphism $N\ra\Pic_{\foX'/S'}$,
thanks to the explicit description of relative Cartier divisors
inside the flat family $\foX'\ra S'$ discussed in the previous section.

The homomorphism $N\ra\Pic_{\foX'/S'}$ is a \emph{marking}
in the sense of Ogus \cite{Ogus 1983}.
He proved in loc.\ cit., Theorem 2.7 that the functor   of isomorphism classes of  $N$-marked K3 surfaces
is representable by  a   algebraic space $S_N$. Note that Ogus works in his paper
under the general assumption $p\geq 3$, which seems appropriate for several arguments
involving quadratic forms.
However, the result concerning the existence of $S_N$
holds true in all characteristics.
Indeed, it is easy to see that the cofibered groupoid $\shF_N$ of $N$-marked families of K3 surfaces
is a stack, and indeed an algebraic stack (= Artin stack). The latter involves Grothendieck's Algebraization Theorem
and openness of versality in the usual way.
According to Rudakov and Shafarevich \cite{Rudakov; Shafarevich 1983}, Section 8, Proposition 3,
the automorphism group of any marked K3 surface is trivial.
By \cite{Ogus 1983}, Lemma 2.2 this carries over to families of marked K3 surfaces.
It follows that the algebraic stack $\shF_N$ is equivalent to its coarse moduli space $S_N$,
which is a nonseparated algebraic space. The upshot is that the algebraic space $S_N$ is
a fine moduli space.

Our $N$-marked flat family $\foX'\ra S'$ induces a morphism $h:S'\ra S_N$.
According to Proposition \ref{generic isomorphism}, this morphism is constant along
all pointed lines, that is, fibers of the canonical projection $S'\ra\PP^1$.
The image $h(S')\subset S_N$ is therefore a 1-dimensional algebraic space.
Any algebraic space is generically a scheme, whence we obtain a rational morphism
$\PP^1\dashrightarrow h(S')$ that   extends as a continuous map on the underlying
topological spaces. A local computation then shows that the rational morphism
extends to a morphism $\PP^1\ra h(S')$ of algebraic spaces.
Pulling back the universal family from $S'$ to $\PP^1$, we obtain the desired family
$\foX\ra\PP^1$ inducing our original family $\foX'\ra S'$.
\qed



\begin{thebibliography}{ccccc}

\bibitem{Artin 1962}
M.~Artin:
Some numerical criteria for contractability of curves on algebraic
surfaces.
Am.\ J.\ Math.\ 84 (1962), 485--496.

\bibitem{Artin 1966}
M.~Artin:
On isolated rational singularities of surfaces.
Am.\ J.\ Math.\ 88 (1966), 129--136.

\bibitem{Artin 1974a}
M.\ Artin:
Supersingular $K3$ surfaces.
Ann.\ Sci.\ \'Ecole Norm.\ Sup.\  7 (1974), 543--567. 

\bibitem{Artin 1974b}
M.\ Artin: 
Algebraic construction of Brieskorn's resolutions. 
J.\ Algebra 29 (1974), 330--348.

\bibitem{Artin 1977}
M.\ Artin:
Coverings of the rational double points in characteristic $p$.
In: W.\ Baily, T.\ Shioda (eds.), 
Complex analysis and algebraic geometry, pp.\ 11--22.
Iwanami Shoten, Tokyo, 1977. 

\bibitem{Bayer; Eisenbud 1995}
D.\ Bayer, D.\ Eisenbud:
Ribbons and their canonical embeddings.
Trans.\ Am.\ Math.\ Soc.\ 347 (1995), 719-756.

\bibitem{Bombieri; Mumford 1976}
E.~Bombieri, D.~Mumford:
Enriques' classification of surfaces in char p.  III.
Invent.\ Math.\ 35  (1976), 197--232.

\bibitem{Bourbaki LIE 1}
N.\ Bourbaki:
Groupes et alg\`ebres de Lie. Chap.\ I: Alg\`ebres de Lie.
Actualit\'es scientifiques et industrielle 1285, Paris, Hermann, 1971.

\bibitem{Bourbaki LIE 4-6}
N.\ Bourbaki:
Groupes et alg\`ebres de Lie. 
Chapitres 4, 5 et 6. 
Masson, Paris, 1981.

\bibitem{Brieskorn 1968}
E.\ Brieskorn:
Die Aufl\"osung der rationalen Singularit\"aten holomorpher Abbildungen.
Math.\ Ann.\  178  (1968) 255--270.

\bibitem{Burns; Rapoport 1975}
D.\ Burns, M.\ Rapoport:
On the Torelli problem for k\"ahlerian K-3 surfaces.
Ann.\ Sci.\ \'Ecole Norm.\ Sup.\  8 (1975),  235--273.

\bibitem{Demazure; Gabriel 1970}
M.\ Demazure, P.\ Gabriel:
Groupes alg\'ebriques.
Masson, Paris, 1970.

\bibitem{Ekedahl 1988}
T.\ Ekedahl:
Canonical models of surfaces of general type in positive characteristic. 
Inst.\ Hautes \'Etudes Sci.\ Publ.\ Math.\ 67 (1988), 97--144.

\bibitem{Giraud 1982}
J.~Giraud:
Improvement of Grauert--Riemenschneider's theorem for a normal surface.
Ann.\ Inst.\ Fourier 32 (1982),  13--23.

\bibitem{Greuel; Kroning 1990}
G.-M.\ Greuel, H.\  Kr\"oning:
Simple singularities in positive characteristic.
Math.\ Z.\ 203 (1990), 339-354.

\bibitem{EGA IVb}
A.\ Grothendieck:
\'El\'ements de g\'eom\'etrie alg\'ebrique IV: \'Etude locale des
sch\'emas et des morphismes de sch\'emas.
Publ.\ Math., Inst.\ Hautes \'Etud.\ Sci.\   24 (1965).

\bibitem{EGA IVd}
A.\ Grothendieck:
\'El\'ements de g\'eom\'etrie alg\'ebrique IV: \'Etude locale des
sch\'emas et des morphismes de sch\'emas.
Publ.\ Math., Inst.\ Hautes \'Etud.\ Sci.\   32 (1967).

\bibitem{SGA 1}
A.~Grothendieck et al.:
Rev\^etements \'etales et groupe fondamental.
Lect.\ Notes Math.\  224,
Springer, Berlin, 1971.

\bibitem{Hartshorne 1994}
R.~Hartshorne:
Generalised divisors on Gorenstein schemes.
K-Theory 8 (1994), 287--339.

\bibitem{Ito 1992}
H.\ Ito:
The Mordell--Weil groups of unirational quasi-elliptic surfaces in characteristic $3$.  
Math.\ Z.\ 211  (1992),  1--39.

\bibitem{Jensen 1978}
S.\ Jensen:
Picard schemes of quotients by finite commutative group schemes.
Math.\ Scand.\  42  (1978), 197--210.

\bibitem{Katsura 1978}
T.\ Katsura:
On Kummer surfaces in characteristic $2$.
In: M.\ Nagata (ed.), Proceedings of the international symposium on 
algebraic geometry, pp.\ 525--542. 
Kinokuniya Book Store, Tokyo, 1978.

\bibitem{Kunz 1986}
E.\ Kunz:
K\"ahler differentials. 
Vieweg, Braunschweig, 1986.

\bibitem{Lipman 1969}
J.\ Lipman:
Rational singularities, with applications to algebraic surfaces and unique factorization.  
Inst.\ Hautes \'Etudes Sci.\ Publ.\ Math.\  36 (1969), 195--279.

\bibitem{Moret-Bailly 1981}
L.\ Moret-Bailly:
Familles de courbes et de varietes abeliennes sur $P^1$.
In: L.\ Szpiro (ed.), S\'eminaire sur les pinceaux de courbes de
genre au moins deux, pp. 109--140.
Ast\'erisque 86 (1981). 

\bibitem{Mumford 1961}
D.~Mumford:
The topology of a normal surface singularity of an algebraic
variety and criterion for simplicity.
Publ.\ Math., Inst.\ Hautes \'Etud.\ Sci.\ 9 (1961), 5--22.

\bibitem{Ogus 1983}
A.\ Ogus:
A crystalline Torelli theorem for supersingular $K3$ surfaces.
In: M.\ Artin, J.\ Tate (eds.),
Arithmetic and geometry, Vol.\ II, pp.\ 361--394.
Progr.\ Math.\ 36.
Birkh\"auser, Boston, 1983.

\bibitem{Oort 1975}
F.\ Oort:
Which abelian surfaces are products of elliptic curves?
Math.\ Ann.\ 214  (1975), 35--47.

\bibitem{Rudakov; Safarevic 1976}
A.\ Rudakov, I.\ Safarevic:
Inseparable morphisms of algebraic surfaces.
Math.\ USSR, Izv.\ 10 (1976), 1205--1237.

\bibitem{Rudakov; Safarevic 1978}
A.\ Rudakov, I.\ Safarevic:
Supersingular $K3$ surfaces over fields of characteristic $2$.
Math.\ USSR, Izv.\ 13 (1979), 147--165.

\bibitem{Rudakov; Shafarevich 1983}
A.\ Rudakov, I.\ Shafarevich:
Surfaces of type $K3$ over fields of finite characteristic. 
In:  I.\ Shafarevich, Collected mathematical papers, pp.\ 657--714.
Springer, Berlin, 1989.

\bibitem{Schroeer 2004}
S.\ Schr\"oer:
Some Calabi--Yau threefolds with obstructed deformations
over the Witt vectors.
Compositio Math.\ 140 (2004), 1579--1592.

\bibitem{Seshadri 1958}
C.\ Seshadri:
Triviality of vector bundles over the affine space $K^2$.
Proc.\ Nat.\ Acad.\ Sci.\ U.S.A.\ 44 (1958), 456--458.

\bibitem{Serre 1965}
J.-P.\ Serre:
Alg\`ebre locale. Multiplicit\'es. 
Lect.\ Notes  Math.\ 11.
Springer, Berlin, 1965.

\bibitem{Serre 1979}
J.-P.\ Serre:
Local fields.
Grad.\ Texts Math.\ 67.
Springer, Berlin, 1979.

\bibitem{Shioda 1978}
T.\ Shioda:
Supersingular $K3$ surfaces.  
In: K.\ Lonsted (ed.), Algebraic geometry, pp.\ 564--591.
Lecture Notes in Math.\ 732. Springer, Berlin, 1979.

\bibitem{Shioda 1974}
T.\ Shioda:
Kummer surfaces in characteristic $2$. 
Proc.\ Japan Acad.\ 50 (1974), 718--722.

\bibitem{Strade; Farnsteiner  1988}
H.\ Strade, R.\  Farnsteiner:
Modular Lie algebras and their representations. 
Monographs and Textbooks in Pure and Applied Mathematics 116. 
Marcel Dekker, New York, 1988.

\bibitem{Wagreich 1970}
P.\ Wagreich:
Elliptic singularities of surfaces.  
Amer.\ J.\ Math.\  92  (1970), 419--454.

\bibitem{Weibel 1994}
C.\ Weibel:
An introduction to homological algebra. 
Cambridge Studies in Advanced Mathematics 38. 
Cambridge University Press, Cambridge, 1994.
\end{thebibliography}
\end{document}